\documentclass[12pt, reqno]{amsproc}

\usepackage[utf8]{inputenc}
\usepackage[T2A]{fontenc}
\usepackage[english]{babel}

\usepackage{amssymb, enumitem}
\usepackage{graphicx}
\usepackage[dvipsnames]{xcolor}
\usepackage{enumitem}
\usepackage[all]{xy}
\usepackage{bm}
\usepackage{bbm}
\usepackage{mathtools}

\usepackage[noadjust]{cite}

\setlist{nosep}

\usepackage[margin=2.9cm]{geometry}

\usepackage{hyperref}
\hypersetup{
    hypertexnames=false,
    colorlinks,
    linkcolor={red!30!black},
    citecolor={blue!50!black},
    urlcolor={blue!80!black}
}

\newtheorem{theorem}[subsection]{Theorem}
\newtheorem{lemma}[subsection]{Lemma}

\newtheorem{corollary}[subsection]{Corollary}

\newtheorem{remark}[subsection]{Remark}

\newtheorem{definition}[subsection]{Definition}

\newtheorem{subtheorem}[subsubsection]{Theorem}
\newtheorem{sublemma}[subsubsection]{Lemma}

\newtheorem{subcorollary}[subsubsection]{Corollary}

\newtheorem{subremark}[subsubsection]{Remark}
\newtheorem{subexample}[subsubsection]{Example}
\newtheorem{subdefinition}[subsubsection]{Definition}

\makeatletter
\@addtoreset{equation}{section}
\@addtoreset{figure}{section}
\@addtoreset{table}{section}
\makeatother

\numberwithin{equation}{section}

\makeatletter
\newcommand\testshape{family=\f@family; series=\f@series; shape=\f@shape.}
\def\myemphInternal#1{\if n\f@shape%
	\begingroup\itshape #1\endgroup\/%
	\else\begingroup\bfseries #1\endgroup%
	\fi}
\def\term{\futurelet\testchar\MaybeOptArgmyemph}
\def\MaybeOptArgmyemph{\ifx[\testchar \let\next\OptArgmyemph
	\else \let\next\NoOptArgmyemph \fi \next}
\def\OptArgmyemph[#1]#2{\index{#1}\myemphInternal{#2}}
\def\NoOptArgmyemph#1{\myemphInternal{#1}}
\makeatother

\newcommand\bN{{\mathbb N}}

\newcommand\bR{{\mathbb R}}

\newcommand\bZ{{\mathbb Z}}

\newcommand\eps{\varepsilon}

\newcommand\id{\mathrm{id}}          
\newcommand\Int{\mathrm{Int}}        
\newcommand\wrm[1]{{\smashoperator{\mathop{\ \wr \ }\limits_{#1}}}} 

\newcommand\BMOD[2]{#1\bmod#2}

\newcommand\AxCrPt{{\text{\rm(H)}}}
\newcommand\AxBd{{\text{\rm(B)}}}

\newcommand\Aman{A}
\newcommand\Bman{B}

\newcommand\Jman{J}
\newcommand\Kman{K}

\newcommand\Mman{M}

\newcommand\Pman{P}
\newcommand\Qman{Q}
\newcommand\Rman{R}

\newcommand\Uman{U}
\newcommand\Vman{V}
\newcommand\Wman{W}
\newcommand\Xman{X}
\newcommand\Yman{Y}
\newcommand\Zman{Z}

\newcommand\Circle{S^1}            

\newcommand\Sphere{S^2}            
\newcommand\Tor[1]{\mathbb{T}^{#1}}

\newcommand\PrjPlane{\bR{P}^2}    

\newcommand\Orb{\mathcal{O}}        
\newcommand\Stab{\mathcal{S}}       
\newcommand\Diff{\mathcal{D}}       
\newcommand\Aut{\mathrm{Aut}}       
\newcommand\Map{\mathrm{Map}}       
\newcommand\FMap{\mathrm{FMap}}       
\newcommand\Maps[2]{\Map(#1,#2)}
\newcommand\MapsFin[2]{\FMap(#1,#2)}

\newcommand\DiffId{\Diff_{\id}}     
\newcommand\StabId{\Stab_{\id}}     

\newcommand\Cinfty{\mathcal{C}^{\infty}}
\newcommand\Cr[1]{\mathcal{C}^{#1}}
\newcommand\Cont[2]{\mathcal{C}\bigl(#1,#2\bigr)}            
\newcommand\Ci[2]{\mathcal{C}^{\infty}(#1,#2)}               
\newcommand\Cid[2]{\mathcal{C}_{\partial}^{\infty}(#1,#2)}   

\newcommand\func{f}

\newcommand\gfunc{g}
\newcommand\dif{h}

\newcommand\Morse{\mathrm{Morse}}

\newcommand\classBC{\classB_{C}}
\newcommand\classBB{\classB_{B}}

\newcommand\Stabilizer[1]{\Stab(#1)}             
\newcommand\StabilizerId[1]{\StabId(#1)}         

\newcommand\Orbit[1]{\Orb(#1)}                   
\newcommand\OrbitPathComp[2]{\Orb_{#2}(#1)}      

\newcommand\AutKRGraphStab[1]{\mathbf{G}}             

\newcommand\flow{\mathbf{F}}

\newcommand\monoArrow{\lhook\joinrel\rightarrow}

\newcommand\xmonoArrow[1]{\lhook\joinrel\xrightarrow{~#1~}}

\newcommand\epiArrow{\rightarrow\!\!\!\!\!\to}

\newcommand\xepiArrow[1]{\xrightarrow{#1}\!\!\!\!\!\to}

\newcommand\restr[2]{#1\vert_{#2}}   

\newcommand\seqC[1]{\mathbf{c}_{#1}}

\newcommand\seqWrm[2]{#1\wr\seqC{#2}}

\newcommand\seqTriv{\seqC{0}} 
\newcommand\pSeq{\mathbf{p}}
\newcommand\qSeq{\mathbf{q}}
\newcommand\aSeq{\mathbf{q}}

\newcommand\kA{K} 
\newcommand\kB{L} 
\newcommand\kC{M} 

\newcommand\OrbfX{\Orbit{\func,\Xman}}        
\newcommand\OrbffX{\OrbitPathComp{\func,\Xman}{\func}}   

\newcommand\FSP[2]{\mathcal{F}(#1,#2)}

\newcommand\Xsp{\mathcal{X}}
\newcommand\Ysp{\mathcal{Y}}

\newcommand\Wsp{\mathcal{W}}

\newcommand\gelg{g}
\newcommand\helg{h}
\newcommand\kelg{k}
\newcommand\lelg{l}
\newcommand\pelg{p}
\newcommand\qelg{q}
\newcommand\xelg{x}
\newcommand\yelg{y}
\newcommand\welg{w}

\newcommand\Ggrp{G}
\newcommand\Hgrp{H}
\newcommand\Kgrp{K}

\newcommand\amapgf{\alpha}

\newcommand\mcirc{\!\circ\!}

\newcommand\qact[3][]{#2#1#3}
\newcommand\qactNoSign[2]{\qact[]{#1}{#2}}
\newcommand\qactDot[2]{\qact[\cdot]{#1}{#2}}

\newcommand\qactBul[2]{\qact[\bullet]{#1}{#2}}
\newcommand\qactCirc[2]{\qact[\circ]{#1}{#2}}
\newcommand\qactSharp[2]{\qact[\,\sharp\,]{#1}{#2}}

\newcommand\gact[2]{\qactNoSign{#1}{#2}}
\newcommand\gpact[2]{\qactCirc{#1}{#2}}
\newcommand\gprod[2]{\qactDot{#1}{#2}}

\newcommand\hact[2]{\qactNoSign{#1}{#2}}
\newcommand\hpact[2]{\qactCirc{#1}{#2}}
\newcommand\hprod[2]{\qactDot{#1}{#2}}

\newcommand\ghact[2]{\qactBul{#1}{#2}}

\newcommand\pconc[2]{\qactSharp{#1}{#2}}
\newcommand\pprod[2]{\qactDot{#1}{#2}}

\newcommand\GwrH{\Ggrp\wr\Hgrp}
\newcommand\Gxorb{\Ggrp\xelg}
\newcommand\Gyorb{\Ggrp\yelg}
\newcommand\Hyorb{\Hgrp\yelg}
\newcommand\GHworb{(\GwrH)\welg}
\newcommand\GKyorb{(\Ggrp/\Kgrp)\yelg}

\newcommand\QXG{\Xman/\Ggrp}
\newcommand\QYH{\Yman/\Hgrp}
\newcommand\QWGH{\Wman/(\GwrH)}

\newcommand\cpathx{\mathbf{\xelg}}

\newcommand\zxmyIso{\eta}
\newcommand\zxmyIsoZ{\zxmyIso_0}

\newcommand\RightShift[1]{R_{#1}}

\newcommand\unit[1]{1_{#1}}

\newcommand\PathSp[3]{\mathcal{P}(#1,#2,#3)}

\newcommand\PXxG{\PathSp{\Xman}{\xelg}{\Gxorb}}

\newcommand\PYyH{\PathSp{\Yman}{\yelg}{\Hyorb}}
\newcommand\PWwGH{\PathSp{\Wman}{\welg}{\GHworb}}
\newcommand\jhom[1]{j_{{}_{#1}}\!}
\newcommand\ehom[1]{\delta_{{}_{#1}}\!}

\newcommand\dhom[1]{\partial_{{}_{#1}}\!}

\newcommand\ehomXx{\ehom{\Xman,\xelg}}
\newcommand\ehomXgx{\ehom{\Xman,\gact{\gelg}{\xelg}}}
\newcommand\ehomXy{\ehom{\Xman,\yelg}}
\newcommand\ehomYy{\ehom{\Yman,\yelg}}

\newcommand\jhomX{\jhom{\Xman}}
\newcommand\ehomX{\ehom{\Xman}}

\newcommand\dhomX{\dhom{\Xman}}

\newcommand\jhomY{\jhom{\Yman}}
\newcommand\ehomY{\ehom{\Yman}}

\newcommand\dhomY{\dhom{\Yman}}

\newcommand\ehomW{\ehom{\Wman}}

\newcommand\xpath{\chi}
\newcommand\ypath{\nu}
\newcommand\wpath{\omega}

\newcommand\classB{\mathcal{B}}

\newcommand\aaa[1]{#1'}
\newcommand\mc[1]{}

\newcommand\agpath{{\mc{Mahogany}\alpha}}
\newcommand\bgpath{{\mc{Mahogany}\beta}}
\newcommand\cgpath{{\mc{Mahogany}\gamma}}

\newcommand\pagpath{{\mc{Mahogany}\aaa{\agpath}}}

\newcommand\ahpath{{\mc{MidnightBlue}\lambda}}
\newcommand\bhpath{{\mc{MidnightBlue}\mu}}
\newcommand\chpath{{\mc{MidnightBlue}\nu}}

\newcommand\aghpath{{\mc{Plum}\psi}}
\newcommand\bghpath{{\mc{Plum}\xi}}
\newcommand\cghpath{{\mc{Plum}\zeta}}

\newcommand\agfunc{{\mc{blue}\phi}}
\newcommand\bgfunc{{\mc{blue}\psi}}

\newcommand\xfunc{{\mc{blue}\agpath}}

\newcommand\agpmap{\aaa{\aaa{\agpath}}}
\newcommand\bgpmap{\aaa{\aaa{\bgpath}}}
\newcommand\cgpmap{\aaa{\aaa{\cgpath}}}
\newcommand\ahpmap{\aaa{\ahpath}}
\newcommand\bhpmap{\aaa{\bhpath}}
\newcommand\chpmap{\aaa{\chpath}}

\newcommand\agmap{\aaa{\agpath}}
\newcommand\bgmap{\aaa{\bgpath}}
\newcommand\cgmap{\aaa{\cgpath}}
\newcommand\ahmap{\ahpath}
\newcommand\bhmap{\bhpath}
\newcommand\chmap{\chpath}

\newcommand\ahm{\mu}
\newcommand\hahm{\delta}

\newcommand\qq[1]{``#1''}

\newcommand\ContU[3]{\mathcal{C}_{#1}(#2,#3)}
\newcommand\ASP{S}
\newcommand\BSP{T}
\newcommand\CSP{X}
\newcommand\apt{s}
\newcommand\bpt{t}
\newcommand\BCSP{\mathcal{U}}
\newcommand\EQU[3]{\pi(#1)}

\newcommand\Dn[1]{D^{#1}}
\newcommand\Sn[1]{S^{#1}}
\newcommand\ptt{q}

\title{Topological actions of wreath products}
\author{Sergiy Maksymenko}
\email{maks@imath.kiev.ua}
\address{Algebra and Topology department, Institute of mathematics, NAS of Ukraine, Str. Tereshchenkivska, 3, Kyiv, 01024, Ukraine}

\begin{document}


\begin{abstract}    
Let $G$ and $H$ be two groups acting on path connected topological spaces $X$ and $Y$ respectively.
Assume that $H$ is finite of order $m$ and the quotient maps $p:X\to X/G$ and $q:Y\to Y/H$ are regular coverings.
Then it is well-known that the wreath product $G\wr H$ naturally acts on $W = X^m\times Y$, so that the quotient map $r:W \to W/(G\wr H)$ is also a regular covering.
We give an explicit description of $\pi_1(W/(G\wr H))$ as a certain wreath product $\pi_1(X/G)\,\wr_{\partial_Y}\pi_1(Y/H)$ corresponding to a \textit{non-effective} action of $\pi_1(Y/H)$ on the set of maps $H\to\pi_1(X/G)$ via the boundary homomorphism $\partial_{Y}:\pi_1(Y/H) \to H$ of the covering map $q$.

Such a statement is known and usually exploited only when $X$ and $Y$ are contractible, in which case $W$ is also contractible, and thus $W/(G\wr H)$ is the classifying space of $G\wr H$.

The applications are given to the computation of the homotopy types of orbits of typical smooth functions $f$ on orientable compact surfaces $M$ with respect to the natural right action of the groups $\mathcal{D}(M)$ of diffeomorphisms of $M$ on $\mathcal{C}^{\infty}(M,\mathbb{R})$.
\end{abstract}

\maketitle

\section{Introduction}\label{sect:intro}
Recall that for two groups $\Ggrp$ and $\Hgrp$ their \term{unrestricted%
\footnote{The \term{restricted wreath product} is defined similarly, but instead of all the set $\Maps{\Hgrp}{\Ggrp}$ one should take its subset $\MapsFin{\Hgrp}{\Ggrp}$ consisting of maps $\alpha\colon \Hgrp\to\Ggrp$ with \qq{finite support} i.e.\ such that $\alpha(\helg)$ differs from the unit of $\Ggrp$ only for finitely many elements $\helg$ of $\Hgrp$.
We will consider only the case when $\Hgrp$ is finite, and in this situation the notions of restricted and non-restricted wreath products coincide.}\!
wreath product} $\Ggrp\wr\Hgrp$ is the semidirect product $\Maps{\Hgrp}{\Ggrp}\rtimes\Hgrp$ corresponding to the natural left action of $\Hgrp$ on the \term{group} $\Maps{\Hgrp}{\Ggrp}$ of all \term{maps} $\Hgrp\to\Ggrp$ (with respect to the point-wise multiplication) by
\[(\helg\alpha)(\gelg) = \alpha(\gelg\helg)\] for $\helg\in\Hgrp$, $\alpha\in\Maps{\Hgrp}{\Ggrp}$, and $\gelg\in\Ggrp$.

It is well known that given a \term{left} action of $\Ggrp$ on a set $\Xman$ and a \term{left} action of $\Hgrp$ on a set $\Yman$, there is a natural left action of the wreath product $\Ggrp \wr \Hgrp$ on $\Wman = \Maps{\Hgrp}{\Xman} \times \Yman$ by $(\alpha,\helg)(\phi,\yelg) = \bigl((\helg\alpha)\phi, \helg\yelg\bigr)$, where $\alpha\in\Maps{\Hgrp}{\Ggrp}$, $\phi\in\Maps{\Hgrp}{\Xman}$, $\helg\in\Hgrp$, $\yelg\in\Yman$, and $(\helg\alpha)\phi\colon \Hgrp\to\Xman$ is given by $\kelg\mapsto\alpha(\kelg\helg)\phi(\helg)$, $\kelg\in\Hgrp$, see Lemma~\ref{lm:wreath_product_action} below.
Usually, such an action is studied for the case when $\Ggrp$ and $\Hgrp$ are finite and \term{freely} act on \term{contractible} topological spaces $\Xman$ and $\Yman$, e.g.~\cite[Theorems~3.3 \& 6.2]{Nakaoka:AM:1961}, \cite[Sect.~2.3.1]{Arnold:PhD:2013}, \cite{Hunton:MPCPS:1990, CastellanaLibman:AM:2009, HopkinsKuhnRavenel:JAMS:2000}.
In that case $\Wman$ can be identified with $\Xman^{m}\times\Yman$, where $m$ is the order of $\Hgrp$, and is contractible as well with respect to the usual product topology.
Moreover, the above action of $\Ggrp \wr \Hgrp$ on $\Wman$ also turns out to be free.
Therefore, the quotient $\Wman/(\Ggrp \wr \Hgrp)$ is usually regarded as the standard model for the classifying space $B(\Ggrp \wr \Hgrp)$ of the group $\Ggrp \wr \Hgrp$.

If $\Hgrp$ is infinite one would need to specify a topology on $\Maps{\Hgrp}{\Xman}$ to make $\Wman$ a suitable topological space.
However, it seems that in general there is no good canonical choice of such a topology, and therefore many authors consider measurable actions and investigate other properties of $\Ggrp \wr \Hgrp$ like amenability, e.g.\ \cite{CornulierStalderValette:TransAMS:2012, BrudeSasyk:CA:2021, BerensteinZamora:NDJFL:2020}.

On the other hand, if $\Xman$ and $\Yman$ are not contractible, and thus $\pi_1(\QXG)$ and $\pi_1(\QYH)$ might differ from $\Ggrp$ and $\Hgrp$, the homotopy type of the quotient $\Wman/(\Ggrp \wr \Hgrp)$ is not well studied.

The main result of the present paper (Theorem~\ref{th:pi1WGwrH}) \term{explicitly} expresses the fundamental group $\pi_1\bigl(\Wman/(\Ggrp \wr \Hgrp)\bigr)$ via the fundamental groups $\pi_1(\QXG)$ and $\pi_1(\QYH)$ as a wreath product corresponding to a certain \term{non-effective action} of $\pi_1(\QYH)$ on the set $\Maps{\Hgrp}{\pi_1(\QXG)}$, under additional assumptions that $\Hgrp$ is finite and each $\Ggrp$ and $\Hgrp$ has a discrete orbit with trivial stabilizer.
The latter conditions are weaker than requiring those actions to be free.
In particular, we get the following statement (Corollary~\ref{cor:coverings}).
\begin{theorem}\sl
Suppose $\Xman$ and $\Yman$ are path connected topological spaces, the actions of $\Ggrp$ and $\Hgrp$ are properly discontinuous, that is the corresponding quotient maps $p\colon\Xman\to\QXG$ and $q\colon\Yman\to\QYH$ are regular coverings, and $\Hgrp$ is finite and consists of $m$ elements.
Let also $\dhomY\colon\pi_1(\QYH)\to\Hgrp$ be the boundary homomorphism of the covering map $q$.
Then $\Ggrp \wr \Hgrp$ freely acts on $\Wman = \Xman^{m}\times\Yman$, the corresponding quotient map $r\colon \Wman\to\QWGH$ is also a regular covering, and we have isomorphisms
\[ \pi_1(\QWGH) \,\cong\, \pi_1(\QXG)\, \wr_{\dhomY} \pi_1(\QYH),\]
see Example~\ref{exmp:wreath_prod_of_homomorphism} for the definition of the right-hand side wreath product, and 
\[
    \pi_k(\QWGH) \cong \pi_k \Wman  \cong  (\pi_k\Xman)^{m} \times \pi_k\Yman, \quad k\geq 2.
\]
\end{theorem}

If $\Xman$ and $\Yman$ are simply connected, then the corresponding boundary homomorphisms $\dhomX\colon \pi_1(\QXG)\cong \Ggrp$ and $\dhomY\colon\pi_1(\QYH)\cong \Hgrp$ of the covering maps $p$ and $q$ are isomorphisms, and we get an isomorphism $\pi_1(\QXG) \wr_{\dhomY} \pi_1(\QYH) \equiv \GwrH$.
This gives a well-known fact that $\pi_1(\QWGH) = \GwrH$.

Applications of this theorem concern with the right action of the group of diffeomorphisms $\Diff(\Mman)$ of a compact surface $\Mman$ on the space $\Ci{\Mman}{\Pman}$ of smooth maps from $\Mman$ to a one-dimensional manifold $\Pman =\bR$ or $\Circle$, see Corollary~\ref{cor:Of_Tn_G}.
We will give one more proof that if $\Mman$ is orientable and differs from $\Sphere$, then there exists a large subset $\FSP{\Mman}{\Pman} \subset \Ci{\Mman}{\Pman}$ (including all Morse functions) such that for every $\func\in\FSP{\Mman}{\Pman}$ its orbit $\OrbitPathComp{\func}{\func}$ has the homotopy type of the quotient $(\Circle)^k/\Ggrp$ of some torus $(\Circle)^k$ by a certain free action of some finite group $\Ggrp$.
In fact, we will present an explicit construction of such an action.
That result was initially established by the author in~\cite{Maksymenko:AGAG:2006, Maksymenko:ProcIM:ENG:2010} for maps $\func\in\FSP{\Mman}{\Pman}$ with trivial $\Ggrp$ (which includes maps taking distinct values at distinct critical points, and in particular generic Morse maps).
Then it was extended by E.~Kudryavtseva~\cite{Kudryavtseva:MathNotes:2012, Kudryavtseva:MatSb:2013, Kudryavtseva:ENG:DAN2016} to the case of non-trivial $\Ggrp$ and maps with singularities locally equivalent to $x^{m}\pm y^2$ (which includes thus all Morse maps), and further by the author in~\cite{Maksymenko:TA:2020} to all $\func\in\FSP{\Mman}{\Pman}$.

\subsection{Structure of the paper}
Section~\ref{sect:prelim} contains preliminary results and definitions.
Section~\ref{sect:exp_law} starts from the exponential law allowing to pass between paths in certain functional spaces and homotopies between maps in those spaces.
Further, for a pointed pair of topological spaces $(\Xman,\Aman,\xelg)$, we discuss the definition of homotopy groups and also homotopy sets $\pi_1(\Xman,\Aman,\xelg)$ and $\pi_0(\Xman,\xelg)$ having in general no ``natural'' groups structures. Section~\ref{sect:pi1XAx} presents several examples when such groups structures exist.
In Section~\ref{sect:WD_actions} we show that if there is an action of some group $\Ggrp$ on $\Xman$ such that the stabilizer of $\xelg$ is trivial and its orbit $\Gxorb$ is discrete (see Definition~\ref{def:WD}), then $\pi_0(\Gxorb)$ can be canonically identified with $\Ggrp$, and $\pi_1(\Xman,\Gxorb,\xelg)$ admits a natural group structure such that the following part of exact sequence of homotopy sets: $\pi_1(\Xman,\xelg)\to\pi_1(\Xman,\Gxorb,\xelg)\to \pi_0(\Gxorb)$ of $(\Xman,\Gxorb,\xelg)$ will consist of homomorphisms, see Lemma~\ref{lm:pi1Xax}.
In particular, many standard facts about the maps of relative $\pi_1$-sets and covering maps can be extended to the case of $\pi_1(\Xman,\Gxorb,\xelg)$.
This construction, of course, is a variant of the concept of \term{monodromy}, see e.g.\ \cite{Zoladek:MG:2006}.

That group structure might be useful as well for discrete dynamical systems.
Namely, if $\phi\colon \Xman\to\Xman$ is a homeomorphism of a path connected topological space $\Xman$, then it yields an action of $\bZ$ on $\Xman$.
Hence, if $\xelg\in\Xman$ is a non-periodic point and its orbit $\bZ\xelg=\{\phi^{k}(\xelg)\}_{\xelg\in\bZ}$ is discrete, then $\pi_1(\Xman,\bZ\xelg,\xelg)$ has a group structure which can be regarded as a certain invariant of $\phi$, see Section~\ref{sect:dyn_syst}.

In Section~\ref{sect:wreath_products} we recall the definition of a wreath product and discuss several particular cases of that construction.
Also, in Section~\ref{sect:short_ex_sequence} we consider certain commutative diagrams related with homomorphisms of wreath products.

Section~\ref{sect:act_wreath_prod} describes a natural action of wreath product $\GwrH$ on $\Wman=\Maps{\Hgrp}{\Xman}\times\Yman$ obtained from actions of groups $\Ggrp$ and $\Hgrp$ on sets $\Xman$ and $\Yman$ respectively, Lemma~\ref{lm:GwrH_W_PD}.
We also prove our main result, Theorem~\ref{th:pi1WGwrH}, which gives an explicit expression of $\pi_1\bigl(\Wman/(\GwrH)\bigr)$ via $\pi_1(\QXG)$ and $\pi_1(\QYH)$, when $\Xman$ and $\Yman$ are topological spaces, and the action of $\Ggrp$ and $\Hgrp$ are \term{weakly discontinuous} at some points $\xelg\in\Xman$ and $\yelg\in\Yman$, see Definition~\ref{def:WD}.

Theorem~\ref{th:pi1WGwrH} is used further in Section~\ref{sect:classes_B} for the proof that a certain class of short exact sequences arises from actions of wreath products on tori, Lemma~\ref{lm:seq_ABC}.
In particular in Remark~\ref{rem:what_to_fix} we also fulfill lacking arguments for~\cite[Theorem 2.5]{Maksymenko:TA:2020}.

\section{Preliminaries}\label{sect:prelim}
In what follows the arrows $\monoArrow$ and $\epiArrow$ will mean a \term{monomorphism} and an \term{epimorphism} respectively.
A \term{diagram} is a functor from the category of partially ordered sets into some category, while a \term{morphism} of diagrams is a natural transformation of the corresponding functors.
For instance, let $\pSeq\colon\Aman_1 \xrightarrow{\alpha_1} \cdots \xrightarrow{\alpha_{k-1}} \Aman_{k}$ and $\qSeq\colon\Bman_1 \xrightarrow{\beta_1}  \cdots \xrightarrow{\beta_{k-1}} \Bman_{k}$ be two sequences of homomorphisms of groups.
Then by a \term{morphism} $\gamma=(\gamma_1,\ldots,\gamma_{k})\colon \pSeq \to \qSeq$ we will mean a collection of homomorphisms $\gamma_i\colon\Aman_i\to \Bman_i$, $i=1,\ldots,k$, making commutative the following diagram:
\[
\aligned
\xymatrix@R=1.1em{
\Aman_1 \ar[r]^-{\alpha_1}  \ar[d]^-{\gamma_1} &
\Aman_2 \ar[r]^-{\alpha_2}  \ar[d]^-{\gamma_2} &
\cdots \ar[r]^-{\alpha_{k-1}} & \Aman_{k} \ar[d]^-{\gamma_k} \\
\Bman_1 \ar[r]^-{\beta_1}   &
\Bman_2 \ar[r]^-{\beta_2}   &
\cdots \ar[r]^-{\beta_{k-1}} & \Bman_{k}
}
\endaligned
\]
A morphism $\gamma$ is an \term{epimorphism} (resp.\ \term{monomorphism}, \term{isomorphism}) if each $\gamma_i$ is so.
We will say that $\gamma$ is a morphism \term{relatively $i$-th term} if $\Aman_i=\Bman_i$ and $\gamma_i=\id_{\Aman_i}$.

\subsection{Homotopies in spaces of maps}\label{sect:exp_law} 
\paragraph{Exponential law}
Let $\BSP,\CSP$ be topological spaces.
By $\Cont{\BSP}{\CSP}$ we will denote the set of all continuous maps from $\BSP$ to $\CSP$.
A \term{compact open topology} on $\Cont{\BSP}{\CSP}$ is a topology generated by sets of the form $[K,V] = \{ \func\in \Cont{\BSP}{\CSP} \mid \func(K)\subset V\}$, where $K$ runs over all compact subsets of $\BSP$ and $V$ runs over all open subsets of $\CSP$.

If $\ASP$ is one more topological space and $F\colon\ASP\times\BSP\to\CSP$ is a continuous map, then for each $\apt\in\ASP$ we define the map $F_{\apt}\colon\BSP\to\CSP$ by $F_{\apt}(\bpt)=F(\apt,\bpt)$, $\bpt\in\BSP$.
The following well known statement is called the \term{exponential law}:
\begin{sublemma}[{e.g.~\cite[Theorem~3.4.8]{Engelking:GT:1989}}]\label{lm:exp_law}\sl
Let $\ASP$, $\BSP$, $\CSP$ be topological spaces such that $\ASP$ is Hausdorff and $\BSP$ is locally compact and Hausdorff.
Then we have a well-defined homeomorphism (called the \term{exponential map}):
\[
    E\colon\Cont{\ASP\times\BSP}{\CSP} \to \Cont{\ASP}{\Cont{\BSP}{\CSP}},
    \qquad
    E(F)(\apt)=F_{\apt},
\]
where $F\in \Cont{\ASP\times\BSP}{\CSP}$, $\apt\in\ASP$, and all the spaces of continuous maps are endowed with the corresponding compact open topologies.
\qed
\end{sublemma}

Given a subset $\BCSP \subset\Cont{\BSP}{\CSP}$, define the following subset of $\Cont{\ASP\times\BSP}{\CSP}$:
\[
    \ContU{\BCSP}{\ASP\times\BSP}{\CSP} :=
    \{
        F\in\Cont{\ASP\times\BSP}{\CSP}
        \mid
        F_{\apt} \in \BCSP \ \text{for all $\apt\in\ASP$}
    \}.
\]
\begin{subcorollary}\label{cor:exp_law}
Under assumptions of Lemma~\ref{lm:exp_law}, for every subset $\BCSP \subset\Cont{\BSP}{\CSP}$, endowed with the induced compact open topology, the map $E$ yields a homeomorphism 
\[ E\colon\ContU{\BCSP}{\ASP\times\BSP}{\CSP} \to \Cont{\ASP}{\BCSP},\]
where $\Cont{\ASP}{\BCSP}$ is regarded as a subset of $\Cont{\ASP}{\Cont{\BSP}{\CSP}}$. 
\end{subcorollary}
\begin{proof}
It suffices to show that $E\bigl(\ContU{\BCSP}{\ASP\times\BSP}{\CSP}\bigr) = \Cont{\ASP}{\BCSP}$.

Let $F\in\ContU{\BCSP}{\ASP\times\BSP}{\CSP}$, so $F_{\apt}\in\BCSP$ for all $\apt\in\ASP$.
Then by Lemma~\ref{lm:exp_law} the map $f=E(F)\colon\ASP\to\Cont{\BSP}{\CSP}$, $f(\apt)=F_{\apt}$, is continuous.
Moreover, $f(\ASP)\subset\BCSP$, by assumption on $F$, i.e.\ $f\in\Cont{\ASP}{\BCSP}$.
Thus, $E\bigl(\ContU{\BCSP}{\ASP\times\BSP}{\CSP}\bigr) \subset \Cont{\ASP}{\BCSP}$.

Conversely, let $f\in\Cont{\ASP}{\BCSP}$.
Then by Lemma~\ref{lm:exp_law}, $f = E(F)$, where $F\in\Cont{\ASP\times\BSP}{\CSP}$ is given by $F(\apt,\bpt)=f(\apt)(\bpt)$.
In particular, $F_{\apt} = f(\apt)\in\BCSP$ for all $\apt\in\ASP$, so $F\in\ContU{\BCSP}{\ASP\times\BSP}{\CSP}$.
Thus, $E\bigl( \ContU{\BCSP}{\ASP\times\BSP}{\CSP} \bigr) = \Cont{\ASP}{\BCSP}$.
\end{proof}

\paragraph{Paths and homotopies}
Let $I=[0;1]$.
Every continuous map $\agpath:I \to \CSP$ is called a \term{path} in $\CSP$.
In this case the points $\agpath(0),\agpath(1)\in\CSP$ are the \term{ends} of $\agpath$, and one also says that $\agpath(0)$ and $\agpath(1)$ are \term{connected by the path $\agpath$}.
If $\agpath(0) = \agpath(1)$, then $\agpath$ is a \term{loop} at the point $\agpath(0)$.
Also, given two subsets $\Aman,\Bman\subset\Xman$, denote by
\[
    \PathSp{\Xman}{\Aman}{\Bman} := \Cont{(I,0,1)}{(\Xman,\Aman,\Bman)}
\]
the set of all paths $\agpath\colon I\to\Xman$ such that $\agpath(0)\in \Aman$ and $\agpath(1)\in\Bman$.

Notice that the relation on $\Xman$ to ``\term{be connected by a path}'' is an equivalence relation.
The corresponding equivalence classes are called \term{path components} of $\Xman$, and the set of all such classes is denoted by $\pi_0\Xman$.

A \term{homotopy} is merely a continuous map $F\colon I\times\BSP\to\CSP$ of a topological product of $\BSP$ by the segment $I$.
In this case, it is also said that $F$ is a \term{homotopy between the maps $F_0,F_1\colon\BSP\to\CSP$}.

Let $\BCSP \subset \Cont{\BSP}{\CSP}$ be a subset.
Then, by definition, each element $F\in\ContU{\BCSP}{I\times\BSP}{\CSP}$ is a homotopy $F\colon I\times\BSP\to\CSP$ satisfying $F_{\apt}\in\BCSP$ for all $\apt\in I$.
We will call such $F$ a \term{homotopy in $\BCSP$}.
Also, say that two maps $\func,\gfunc\in\BCSP$ are \term{homotopic in $\BCSP$} if there exists a homotopy $F$ in $\BCSP$ such that $F_0=\func$ and $F_1=\gfunc$.
Evidently, the relation to ``\term{be homotopic in $\BCSP$}'' is an equivalence relation on $\BCSP$, and we will denote by $\EQU{\BCSP}{\BSP}{\CSP}$ the corresponding set of equivalence classes, called \term{homotopy classes}.

\begin{subcorollary}
Suppose $\BSP$ is locally compact and Hausdorff.
Then for every subset $\BCSP \subset \Cont{\BSP}{\CSP}$ the exponential map $E$ induces a bijection 
\begin{equation}\label{equ:Pi0U_pi0U}
    \EQU{\BCSP}{\BSP}{\CSP} \equiv \pi_0 \BCSP
\end{equation}
between the set homotopy classes of maps in $\BCSP$ and the set of path components of $\BCSP$ endowed with the compact open topology.
\end{subcorollary}
\begin{proof}
By Corollary~\ref{cor:exp_law}, we have a homeomorphism $E\colon\ContU{\BCSP}{I\times\BSP}{\CSP} \cong\Cont{I}{\BCSP}$ of the \term{space of homotopies in $\BCSP$} onto \term{the space of paths in $\BCSP$}.
In particular, two maps $\func,\gfunc\in\BCSP$ are homotopic in $\BCSP$ if and only if they are conected by a path in $\BCSP$, i.e.\ belong to the same path component of $\BCSP$.
This gives the required bijection~\eqref{equ:Pi0U_pi0U}. 
\end{proof}

Let $k\geq1$, $\BSP_1, \ldots, \BSP_k\subset\BSP$ and $\CSP_1, \ldots, \CSP_k\subset\CSP$ be two collections of subsets, and $\BCSP:= \Cont{(\BSP, \BSP_{k},\ldots,\BSP_1)}{(\CSP, \CSP_{k},\ldots,\CSP_1)}$ be the subspace of $\Cont{\BSP}{\CSP}$ consisting of maps $\gamma\colon\BSP\to\CSP$ of $(k+1)$-tuples, i.e.\ $\gamma(\BSP_i)\subset\CSP_i$ for all $i=1,\ldots,k$.
Then the corresponding set of homotopy classes $\EQU{\BCSP}{\BSP}{\CSP}$ is usually denoted by
$[(\BSP, \BSP_{k},\ldots,\BSP_1), (\CSP, \CSP_{k},\ldots,\CSP_1)]$.

\paragraph{Homotopy groups and sets}
For the convenience of the reader not familiar with homotopy groups we will briefly recall their definition, see e.g.\ \cite[Chapter~4.1]{Hatcher:AlgTop:2002}, and further discuss the cases when they are not actually groups.

For $n\geq1$ let $\Dn{n}$ be the unit $n$-disk in $\bR^{n}$ centered at the origin, $\Sn{n-1} = \partial\Dn{n}$ be the corresponding $n$-sphere, and $\ptt=(1,0,\ldots,0)\in\Sn{n-1}$.
Identifying $\bR^{n}$ with $\bR^{n}\times0\subset\bR^{n+1}$ we can assume that $\ptt$ is a common point of all those spheres and disks, and therefore we will use the same point $\ptt$ for all $n$.

Now let $\Xman$ be a topological space, $\Aman\subset\Xman$ a subset, and $\xelg \in \Aman$ a point.
Then the following sets of homotopy classes:
\begin{align}
    \label{equ:pinXx}\pi_{n}(\Xman,\xelg)       &:= [(\Sn{n},\ptt), (\Xman,\xelg)],                \ n\geq0, \\
    \label{equ:pinXAx}\pi_{n}(\Xman,\Aman,\xelg) &:= [(\Dn{n},\Sn{n-1},\ptt), (\Xman,\Aman,\xelg)], \ n\geq1,
\end{align}
are called the \term{$n$-th homotopy set} of $\Xman$ and $(\Xman,\Aman)$ respectively at the point $\xelg$.
Notice that there is an infinite (to the left) sequence of maps called the \term{long exact sequence of $(\Xman,\Aman,\xelg)$}, e.g.~\cite[Theorem~4.3]{Hatcher:AlgTop:2002}:
\begin{multline*}
    \cdots \to
    \pi_{n}(\Aman,\xelg)
    \xrightarrow{i_{n}}
    \pi_{n}(\Xman,\xelg)
    \xrightarrow{j_{n}}
    \pi_{n}(\Xman,\Aman,\xelg)
    \xrightarrow{\partial_n}
    \pi_{n-1}(\Aman,\xelg)
    \to \cdots  \\
    \cdots \to
    \pi_{1}(\Aman,\xelg)       \xrightarrow{i_1}
    \pi_{1}(\Xman,\xelg)       \xrightarrow{j_1}
    \pi_{1}(\Xman,\Aman,\xelg) \xrightarrow{\partial_1}
    \pi_{0}(\Aman,\xelg)       \xrightarrow{i_0}
    \pi_{0}(\Xman,\xelg)
\end{multline*}
and defined as follows.
Let $i:\Aman\subset\Xman$ be a natural inclusion.
\begin{itemize}[leftmargin=*]
\item
Then $i_{n}$ associates to the homotopy class of each map $\gamma\colon(\Sn{n},\ptt) \to (\Aman,\xelg)$ the homotopy class of the composition $i\circ\gamma\colon(\Sn{n},\ptt) \xrightarrow{\gamma} (\Aman,\xelg) \xrightarrow{i} (\Xman,\xelg)$.

\item
Also, note that there is a continuous map $\phi\colon\Dn{n}\to\Sn{n}$ such that $\phi(\partial\Dn{n})=\ptt$ and $\phi$ homemorphically maps $\Int{\Dn{n}}$ onto $\Sn{n}\setminus\ptt$.
In particular, $\phi$ can be regarded as a map of triples $\phi\colon(\Dn{n},\Sn{n-1},\ptt)\to(\Sn{n},\ptt,\ptt)$.
Then $j_n$ associates to the homotopy class of each map of pairs $\gamma\colon(\Sn{n},\ptt) \to (\Xman,\xelg)$ the homotopy class of the following composition regarded as a map of triples:
\[
    \gamma\circ\phi\colon
    (\Dn{n},\Sn{n-1},\ptt)
        \xrightarrow{\phi}
    (\Sn{n},\ptt,\ptt)
        \xrightarrow{\gamma}
    (\Xman,\Aman,\xelg).
\]

\item
Finally, $\partial_n$ associates to the homotopy class of the map of triples $\gamma\colon (\Dn{n},\Sn{n-1},\ptt) \to (\Xman,\Aman,\xelg)$ the homotopy class of its restriction $\restr{\gamma}{(\Sn{n-1},\ptt)}\colon(\Sn{n-1},\ptt)\to(\Aman,\xelg)$ regarded as a map of pairs.
\end{itemize}
It is known that the above long sequence is \term{exact} in the ``homotopy'' sense, i.e.\ if $\Pman\xrightarrow{a}\Qman\xrightarrow{b}\Rman$ is a part of that sequence and $\gamma\in\Qman$ is any element, then $b(\gamma)$ is homotopic in $\Rman$ to a constant map if and only if $\gamma=a(\beta)$ for some $\beta\in\Pman$.

Moreover, the sets $\pi_{n}(\Aman,\xelg)$, $\pi_{n}(\Xman,\xelg)$ for $n\geq1$ and $\pi_{n}(\Xman,\Aman,\xelg)$ for $n\geq2$ have natural group structures such that the maps $i_n$ for $n\geq1$ and $j_n,\partial_n$ for $n\geq2$ are homomorphisms and the above infinite sequence is also \term{exact} in the ``algebraic'' sense up to the term $\pi_1(\Xman,\xelg)$, i.e.\ if $\Pman\xrightarrow{a}\Qman\xrightarrow{b}\Rman$ is a part of that sequence, and $\Qman$ stands before $\pi_1(\Xman,\xelg)$ in that sequence, then $\mathrm{image}(a)=\ker(b)$.

The reader is referred to~\cite[Chapter~4.1]{Hatcher:AlgTop:2002} for the definition of multiplication in those homotopy sets, and we will recall now only the multiplication in $\pi_1(\Xman,\xelg)$ and discuss other homotopy sets of dimension $1$ and $0$.
It will be convenient to replace in the definitions $\Dn{1}=[-1;1]$ with $I=[0;1]$.
Then $\Sn{0}:=\partial I = \{0,1\}$.

\paragraph{Fundamental group $\pi_1(\Xman,\xelg)$}
Define the following map $\phi\colon I\to\Circle$, $\phi(t)=(\cos(2\pi t), \sin(2\pi t))$.
Evidently, $\phi(0)=\phi(1) = \ptt$, and $\phi$ homeomorphically maps $(0;1)$ onto $\Sn{1}\setminus\ptt$.
One easily checks that we have a bijection
\begin{gather*}
    \phi_{*}\colon \Cont{(\Sn{1},\ptt)}{(\Xman,\xelg)} \to \Cont{(I,\partial I)}{(\Xman,\xelg)} \equiv \PathSp{\Xman}{\xelg}{\xelg},
    \\
    \phi_{*}(\gamma) = \gamma\circ\phi\colon I \xrightarrow{\phi} \Sn{1} \xrightarrow{\gamma} \Xman,
\end{gather*}
which also induces a bijection on the corresponding homotopy classes of maps.
Hence, we get another description of $\pi_1(\Xman,\xelg)$ as the set of homotopy classes of \term{loops at $\xelg$}:
\[ 
    \pi_{1}(\Xman,\xelg) 
    \stackrel{\eqref{equ:pinXx}}{:=} 
    [(\Sn{1},\ptt), (\Xman,\xelg)] 
    \stackrel{\phi_{*}}{\equiv} 
    [(I,\partial I), (\Xman,\xelg)].
\]

Say that two paths $\agpath,\bgpath\colon I \to\Xman$ are \term{composable}, if $\agpath(1) = \bgpath(0)$.
Then for a pair of composable paths $\agpath,\bgpath\colon I \to\Xman$ one can define their \term{composition} $\pconc{\agpath}{\bgpath}\colon I\to\Xman$ by the following standard formula:
\begin{equation}\label{equ:composable_paths}
(\pconc{\agpath}{\bgpath})(t) =
\begin{cases}
\agpath(2t), & t\in[0;\tfrac{1}{2}], \\
\bgpath(2t-1), & t\in[\tfrac{1}{2};1].
\end{cases}
\end{equation}
In particular, if $\alpha, \beta \in \PathSp{\Xman}{\xelg}{\xelg}$ are \term{loops} at $\xelg$, then $\pconc{\agpath}{\bgpath}\in\PathSp{\Xman}{\xelg}{\xelg}$ is also a loop at $\xelg$.
It is well known and is easy to see that the homotopy class of $\pconc{\agpath}{\bgpath}$ in $\PathSp{\Xman}{\xelg}{\xelg}$ depends only on the homotopy classes of $\alpha$ and $\beta$, and the operation of composition of loops induces a certain \term{group} operation on the set $\pi_1(\Xman,\xelg)$ of their homotopy classes.
The group $\pi_1(\Xman,\xelg)$ is called the \term{fundamental group} of $\Xman$ at $\xelg$.

However, in general, neither of the latter three terms of the long exact sequence:
\begin{equation}\label{equ:long_ex-seq__end}
    \cdots \xrightarrow{j_1} \pi_{1}(\Xman,\Aman,\xelg) \xrightarrow{\partial_1}\pi_{0}(\Aman,\xelg)  \xrightarrow{i_0}\pi_{0}(\Xman,\xelg)
\end{equation}
has some natural group structure under which the corresponding arrows are homomorphisms.

\paragraph{The set $\pi_1(\Xman,\Aman,\xelg)$}
Evidently, for $n=1$ the general definition~\eqref{equ:pinXAx} reduces to the following one: $\pi_1(\Xman,\Aman,\xelg) = [(I,1,0),(\Xman,\Aman,\xelg)]$.
Note that, if $\agpath,\bgpath\colon(I,1,0)\to(\Xman,\Aman,\xelg)$ are two paths started at $\xelg$ and finished in $\Aman$, then in general, $\agpath$ and $\bgpath$ are not composable.
This illustrates an absence of some ``evident'' group structure on $\pi_1(\Xman,\Aman,\xelg)$ related with compositions of paths.

\paragraph{The set $\pi_0(\Xman,\xelg)$}
Due to definition~\eqref{equ:pinXx}, $\pi_0(\Xman,\xelg) := [(\Sn{0},\ptt),(\Xman,\xelg)]$.
Since every continuous map $\gamma\colon(\Sn{0}, 0)\to(\Xman,\xelg)$ sends $0$ to $\xelg$, it is uniquely determined by its value $\gamma(1)\in\Xman$.
Moreover, the correspondence $\Cont{(\Sn{0}, 0)}{(\Xman,\xelg)}\to\Xman$, $\gamma\mapsto\gamma(1)$, is evidently a bijection.
Also, if $F\colon I \times \Sn{0}\to \Xman$ is a homotopy in $\Cont{(\Sn{0}, 0)}{(\Xman,\xelg)}$, then $F(I\times 0) = \{\xelg\}$, and $\restr{F}{I\times 1}$ is just a path between two points $F_0(1)$ and $F_1(1)$.
This implies that $\pi_0(\Xman,\xelg)$ can be identified with the set of path components of $\Xman$.

From that point of view the presence of a base point $\xelg$ in $\pi_0(\Xman,\xelg)$ might look artificial.
Therefore, sometimes it is omitted and this agrees with the notation $\pi_0\Xman$ from Section~\ref{sect:exp_law}.
However, for the ``exactness'' of the above long exact sequence of $(\Xman,\Aman,\xelg)$ the base point is essential.

Let us also describe the maps from~\eqref{equ:long_ex-seq__end}.
\begin{itemize}[leftmargin=*, itemsep=1ex]
\item
The mapping $j_1\colon\pi_1(\Xman,\xelg)\to \pi_{1}(\Xman,\Aman,\xelg)$ associates to each loop $\alpha\colon(I,\partial I)\to(\Xman,\ptt)$ at $\xelg$ the homotopy class of $\alpha$ regarded as a path $\alpha\colon(I,1,0)\to(\Xman,\Aman,\xelg)$.

\item
The map $\partial_1\colon\pi_{1}(\Xman,\Aman,\xelg)\to\pi_0(\Aman,\xelg)$ associates to each homotopy class of a path $\gamma:(I,1,0)\to(\Xman,\Aman,\xelg)$ the path component of $\Aman$ of the point $\gamma(1)\in\Aman$.

\item
Finally, the map $i_0\colon\pi_{0}(\Aman,\xelg)\to\pi_0(\Xman,\xelg)$ associates to each path component $\Aman'$ of $\Aman$ the path component of $\Xman$ containing $\Aman'$.
\end{itemize}

\begin{subcorollary}\label{cor:pi0_PathsAndLoops}
For a subset $\Aman\subset\Xman$ and a point $\xelg\in\Xman$ there are the following identifications:
\begin{align*}
    \pi_0\PathSp{\Xman}{\xelg}{\Aman} &= \pi_1(\Xman,\Aman,\xelg), &
    \pi_0\Omega(\Xman,\xelg)          &= \pi_1(\Xman,\xelg),
\end{align*}
where $\Omega(\Xman,\xelg) := \PathSp{\Xman}{\xelg}{\xelg}$ is the space of loops at $\xelg$.
\end{subcorollary}
\begin{proof}
The proof is a direct consequence of definitions and~\eqref{equ:Pi0U_pi0U}:
\[
    \pi_1(\Xman,\Aman,\xelg) :=
    [(I,1,0),(\Xman,\Aman,\xelg)]
    \equiv
    \EQU{\PathSp{\Xman}{\xelg}{\Aman}}{I}{\CSP}
    \stackrel{\eqref{equ:Pi0U_pi0U}}{=}
    \pi_0\PathSp{\Xman}{\xelg}{\Aman}.
\]
Similarly,
$    \pi_0(\Xman,\xelg) :=
    [(I,1,0),(\Xman,\xelg,\xelg)]
    \equiv
    \EQU{\Omega(\Xman,\xelg)}{I}{\CSP}
    \stackrel{\eqref{equ:Pi0U_pi0U}}{=}
    \pi_0\Omega(\Xman,\xelg).
$
\end{proof}

\subsection{Group structure on $\pi_1(\Xman,\Aman,\xelg)$}\label{sect:pi1XAx}

We will consider here several situations in which $\pi_1(\Xman,\Aman,\xelg)$, $\pi_0(\Xman,\xelg)$,  $\pi_0(\Aman,\xelg)$ are groups and the corresponding arrows between them are homomorphisms.

\begin{subexample}\rm
Let $\Xman$ be a topological group.
Then the path component $\Xman_{e}$ of the unit $e$ of $\Xman$ is a normal subgroup of $\Xman$, and one can naturally identify $\pi_0\Xman$ with the quotient group $\Xman/\Xman_{e}$.
One can also define a point-wise multiplication of loops $\alpha,\beta\colon(I,\partial I)\to(\Xman,e)$, and it is well known and is easy to see that on the level of homotopy classes that multiplication coincides with the multiplication in $\pi_1(\Xman,e)$.
More generally, if $\Aman\subset\Xman$ is a subgroup, then one can also define a point-wise multiplication of paths $\alpha,\beta\colon(I,0,1)\to(\Xman,e,\Aman)$ and on the level of homotopy classes this turns $\pi_1(\Xman,\Aman, e)=\pi_0\PathSp{\Xman}{e}{\Aman}$ into a group such that the corresponding exact sequence of $(\Xman,\Aman,e)$:
\[
\cdots \to
\pi_1(\Aman,e)       \xrightarrow{i_1}
\pi_1(\Xman,e)       \xrightarrow{j_1}
\pi_1(\Xman,\Aman,e) \xrightarrow{\partial_1}
\pi_0(\Aman,e)       \xrightarrow{i_0}
\pi_0(\Xman,e)
\]
consists of homomorphisms.
It is well known and easy that $\pi_1(\Xman,e)$ and $\pi_1(\Aman,e)$ are abelian, and the image of $j_1$ is contained in the center of $\pi_1(\Xman,\Aman,e)$.
\end{subexample}

\begin{subexample}\rm
Let $p\colon\Xman\to\Yman$ be a fibration between path connected spaces, i.e.\ $p$ satisfies homotopy lifting property.
Let also $\xelg\in\Xman$, $\yelg=p(\xelg)\in\Yman$, and $\Aman = p^{-1}(\yelg)$.
Then that homotopy lifting property implies that $p$ induces a bijection $p_{k}\colon\pi_{k}(\Xman,\Aman,\xelg) \cong \pi_{k}(\Yman,\yelg)$ being an isomorphism of groups for $k\geq2$.
For $k=1$ we have a bijection $p_1$ between $\pi_{1}(\Xman,\Aman,\xelg)$ and the \term{group} $\pi_{1}(\Yman,\yelg)$.
This allows to endow $\pi_1(\Xman,\Aman,\xelg)$ with a groups structure from $\pi_{1}(\Yman,\yelg)$ via $p_1$.
However, the boundary map $\partial_1\colon\pi_1(\Xman,\Aman,\xelg)\to \pi_0(\Aman,\xelg)$ is still not a homomorphism, since $\pi_0(\Aman,\xelg)$ is not a group.
\end{subexample}

\begin{subexample}\rm
Suppose that in the previous example $p\colon\Xman\to\Yman$ is a regular covering map, which means that $\Aman$ is discrete, so $p_1\colon\pi_{1}(\Xman,\xelg) \to \pi_{1}(\Yman,\yelg)$ is injective, and the image of $p_1$ is a normal subgroup of $\pi_{1}(\Yman,\yelg)$.
Then the quotient group $\Ggrp = \pi_{1}(\Yman,\yelg)/p_1\bigl(\pi_{1}(\Xman,\xelg)\bigr)$ can be naturally identified with $\Aman=\pi_0\Aman$ so that the boundary map $\partial_1\colon\pi_1(\Xman,\Aman,\xelg) \to \pi_0(\Aman,\xelg) \equiv \Ggrp$ becomes a homomorphism.

In particular, the non-trivial part $\pi_1(\Xman,\xelg) \xmonoArrow{j_1} \pi_1(\Xman,\Aman,\xelg) \xepiArrow{\partial_1} \pi_0(\Aman,\xelg)$ of the long exact sequence of pair $(\Xman,\Aman,\xelg)$ turns into $\pi_1(\Xman,\xelg) \xmonoArrow{j_1} \pi_1(\Yman,\yelg) \xepiArrow{\partial_1} \Ggrp$.
In this case we have a natural action of $\Ggrp$ on $\Xman$ so that $\Xman/\Ggrp$ can be identified with $\Yman$ and $\Aman$ with the orbit $\Gxorb$ of $\xelg$.
\end{subexample}

\subsection{Weakly discontinuous actions}\label{sect:WD_actions}
Suppose now that a discrete group $\Ggrp$ acts from the left on a topological space $\Xman$ by homeomorphisms.
For a point $\xelg\in\Xman$ its $\Ggrp$-orbit will be denoted by $\Gxorb$.
Our aim is to show that the standard arguments from the theory of covering spaces allow to prove that $\pi_1(\Xman,\Gxorb,\xelg)$ still has a group structure under more general settings.
This should probably be known for specialists, however the author did not find any such exposition in the literature.

For a point $a\in\Xman$ denote by $[a]$ its path component in $\Xman$.

\begin{subdefinition}\label{def:WD}\sl
Say that a $\Ggrp$-action is \term{weakly discontinuous at a point $\xelg\in\Xman$} (WD at $\xelg$) if it satisfies either of the following equivalent properties:
\begin{enumerate}[leftmargin=*, itemsep=0.5ex, label={\rm(WD\arabic*)}]
\item
the natural map $\sigma_{x}\colon\Ggrp \xrightarrow{~\gelg\,\mapsto\,\gact{\gelg}{\xelg}~} \Gxorb \xrightarrow{~a\,\mapsto\,[a]~} \pi_0(\Gxorb,\xelg)$ is a bijection;
\item
the path component of $\xelg$ in $\Gxorb$ consists of that point $\xelg$ only, and the stabilizer of $\xelg$ is trivial;
\item
every continuous path $\alpha\colon I\to\Gxorb$ is constant, and if $\gact{\gelg}{\xelg}=\gact{\gelg'}{\xelg}$ for some $\gelg,\gelg'\in\Ggrp$, then $\gelg=\gelg'$.
\end{enumerate}
\end{subdefinition}
One also checks that every WD-action at $\xelg$ is effective (since distinct elements of $\Ggrp$ differently act on $\Gxorb$) and is also WD at any other point $\yelg$ of the orbit $\Gxorb$.

\begin{subdefinition}\label{def:PD}\sl
Say that $\Ggrp$ acts \term{properly discontinuous} (PD) if either of the following equivalent conditions holds:
\begin{enumerate}[leftmargin=*, label={\rm(PD\arabic*)}]
\item\label{enum:def:PD:covering}
the action is free and the quotient map $p\colon\Xman\to\QXG$ is a covering, where $\QXG$ is endowed with the corresponding quotient topology;
\item\label{enum:def:PD:wander_nbh}
each $\xelg\in\Xman$ has an open neighborhood $\Vman_{\xelg}$ such that $\gact{\gelg}{\Vman_{\xelg}} \cap \gact{\helg}{\Vman_{\xelg}}=\varnothing$ for $\gelg\not=\helg$.
\end{enumerate}
\end{subdefinition}
A neighborhood $\Vman_{\xelg}$ in~\ref{enum:def:PD:wander_nbh} is called \term{wandering (with respect to this action of $\Ggrp$)} or simply \term{$\Ggrp$-wandering}.

Evidently, every PD action is WD at each point $\xelg\in\Xman$ and the restriction of $p$ to every $\Ggrp$-wandering neighborhood is a homeomorphism onto.
Moreover, it is well known and is easy to see that if $\Xman$ is Hausdorff and $\Ggrp$ is a finite group freely acting on $\Xman$, then this action is also PD, e.g.\ \cite[11.1.3]{Brown:TG:2006}.

In fact, there are several definitions of \term{properly} discontinuous actions involving compact subspaces and adopted to the actions on locally compact Hausdorff spaces, e.g.\ \cite[I, page 4, Remark]{Charlap:BG:1986}.
We will use only the one equivalent to~\ref{enum:def:PD:covering}.

\subsection{Group structure on $\pi_0\PXxG$}
Let $\Xman$ be a path connected topological space and $\xelg\in\Xman$.
Suppose we are given an action of a group $\Ggrp$ on $\Xman$ being WD at $\xelg$.

1) Then there is a natural map
\[
    \ehomXx\colon\PXxG \to \Ggrp
\]
defined as follows.
Let $\agpath\in\PXxG$, so it is a path $\alpha\colon I \to \Xman$ such that $\agpath(0)=\xelg$ and $\agpath(1) = \gact{\gelg}{\xelg} \in \Gxorb$ for some $\gelg\in\Ggrp$.
Since the action is WD at $\xelg$ the stabilizer of $\xelg$ is trivial, and therefore such $\gelg$ is unique.
Moreover, it also depends only on the homotopy class $[\agpath]$ of $\agpath$ in $\pi_0\PXxG\equiv\pi_1(\Xman,\Gxorb,\xelg)$.

Indeed, if $\{\agpath_t\}_{t\in[0;1]} \subset\PXxG$ is a homotopy in $\PXxG$, then $\{\agpath_t(1) \}_{t\in[0;1]}$ is a path in $\Gxorb$ which must be constant, since the path components of $\Ggrp$ are singletons (by WD property at $\xelg$).
Hence, the correspondence $\agpath \mapsto \gelg$ is a well-defined map $\ehomXx\colon\pi_1(\Xman,\Gxorb,\xelg)\to\Ggrp$ such that $\agpath(1) = \gact{\ehomXx([\alpha])}{\xelg}$.

2) Further, one can define the following operation of concatenation of paths:
\[
    \PXxG \times \PXxG \to \PXxG
\]
in the following way.
Let $\agpath,\bgpath\in\PXxG$ and $\gelg = \ehomXx(\agpath)$, so $\agpath(1)=\gact{\gelg}{\xelg}$ for a unique $\gelg\in\Ggrp$ which can be regarded as a homeomorphism of $\Xman$.
Then the composition $\gpact{\gelg}{\bgpath}\colon I \xrightarrow{\bgpath} \Xman \xrightarrow{\gelg} \Xman$ is also a path in $\Xman$, and $(\gpact{\gelg}{\bgpath})(0)=\gact{\gelg}{\xelg}=\agpath(1)$.
Hence, $\agpath$ and $\gpact{\gelg}{\bgpath}$ are composable, and we define the product of elements $\agpath,\bgpath\in\PXxG$ by
\begin{equation}\label{equ:prod_xpath_xpath1}
    \pprod{\agpath}{\bgpath} :=
    \pconc{\agpath}{(\gpact{\gelg}{\bgpath})} \equiv
    \pconc{\agpath}{(\gpact{\ehomXx(\agpath)}{\bgpath})}.
\end{equation}

\begin{sublemma}\label{lm:pi1Xax}\sl
Suppose $\Xman$ is path connected and the action of $\Ggrp$ on $\Xman$ is WD at some $\xelg\in\Xman$.
Then the following statements hold.
\begin{enumerate}[leftmargin=*, label={\rm(\arabic*)}, topsep=0.5ex, itemsep=0.5ex]
\item\label{enum:lm:pi1Xax:ehomX}
$\pi_1(\Xman,\Gxorb,\xelg)$ has a group structure such that the composition
\[
    \ehomXx = \sigma_x^{-1}\circ\partial_1:
        \pi_1(\Xman,\Gxorb,\xelg) \xrightarrow{~\partial_1~}
        \pi_0(\Gxorb,\xelg)       \xrightarrow{~\sigma_x^{-1}~}
        \Ggrp
\]
is a homomorphism.
In particular, we get a short exact sequence:
\begin{equation}\label{equ:short_ex_seq_X_Gx_x}
   \pi_1(\Xman,\xelg)        \xmonoArrow{~j_1~}
   \pi_1(\Xman,\Gxorb,\xelg) \xrightarrow{~\ehomXx~}
   \Ggrp.
\end{equation}

\item\label{enum:lm:pi1Xax:iso_by_path}
Suppose that this $\Ggrp$-action is also WD at some other point $\yelg\in\Xman$.
Let $\gamma\colon I\to\Xman$ be any path with $\gamma(0)=\xelg$ and $\gamma(1)=\yelg$.
Then the natural bijection
\[
    \gamma_{*}\colon\pi_1(\Xman,\Gxorb,\xelg) \to\pi_1(\Xman,\Gyorb,\yelg),
    \qquad
    \gamma_{*}(\agpath)=[\pconc{\gamma^{-1}}{\pconc{\agpath}{(\ehomXx(\agpath)\circ\gamma})}],
\]
is an isomorphism of groups inducing an isomorphism of the following short exact sequences:
\[
\xymatrix{
\ \pi_1(\Xman,\xelg)          \ \ar[d]_-{\cong}^{\gamma_{*}}    \ar@{^(->}[r] &
\ \pi_1(\Xman,\Gxorb,\xelg)   \ \ar[d]^-{\gamma_{*}}_-{\cong} \ar@{->>}[r]^-{\ehomXx} &
\ \Ggrp                       \ \ar@{=}[d] \\
\ \pi_1(\Xman,\yelg)          \ \ar@{^(->}[r] &
\ \pi_1(\Xman,\Gyorb,\yelg)   \ \ar@{->>}[r]^-{\ehomXy} &
\ \Ggrp                       \
}
\]

\item\label{enum:lm:pi1Xax:equivariant_maps}
Assume that another group $\Hgrp$ acts on a path connected topological space $\Yman$ and that action is WD at some $\yelg\in\Yman$.
Let also $\func\colon\Xman\to\Yman$ be a continuous map such that $\func(\xelg)=\yelg$.
Suppose also that there exists a homomorphism $\phi\colon \Ggrp\to\Hgrp$ such that $\func$ is \term{$\phi$-equivariant on the orbit $\Gxorb$ of $\xelg$}, i.e.\ $\func(\gact{\gelg}{\xelg})=\hact{\phi(\gelg)}{\yelg}$ for all $\gelg\in\Ggrp$.
In particular, $\func(\Gxorb) \subset \Hyorb$.
Then the natural map
\begin{equation}\label{equ:f_pi1XGx_pi1YHy}
    \func_{*}\colon\pi_1(\Xman,\Gxorb,\xelg) \to \pi_1(\Yman,\Hyorb,\yelg),
    \qquad
    \func_{*}([\agpath]) = [\func\circ\agpath],
\end{equation}
\term{is not necessarily a homomorphism} however it induces the following commutative diagram in which the left vertical arrow $\func_{*}$ is a homomorphism:
\begin{equation}\label{equ:func_morphism_of_sequences}
\begin{aligned}
    \xymatrix{
    \ \pi_1(\Xman,\xelg)          \ \ar[d]^-{\func_{*}}  \ar@{^(->}[r] &
    \ \pi_1(\Xman,\Gxorb,\xelg)   \ \ar[d]^-{\func_{*}} \ar@{->>}[r]^-{\ehomXx} &
    \ \Ggrp                       \ \ar[d]^-{\phi} \\
    \ \pi_1(\Yman,\yelg)          \ \ar@{^(->}[r] &
    \ \pi_1(\Xman,\Hyorb,\yelg)   \ \ar@{->>}[r]^-{\ehomYy} &
    \ \Hgrp                       \
    }
\end{aligned}
\end{equation}

If $\func$ is $\phi$-equivariant on all of $\Xman$, i.e.\ $\func(\gact{\gelg}{\xelg'})=\hact{\phi(\gelg)}{\func(\xelg')}$ for all $\gelg\in\Ggrp$ and $\xelg'\in\Xman$, then~\eqref{equ:f_pi1XGx_pi1YHy} is a homomorphism, so~\eqref{equ:func_morphism_of_sequences} is a morphism of short exact sequences, i.e.\ all vertical arrows there are homomorphisms.
\end{enumerate}
\end{sublemma}
\begin{proof}
\ref{enum:lm:pi1Xax:ehomX}
A standard verification shows that the above operation $\sharp$ at the level of homotopy classes of paths in $\PXxG$ turns $\pi_0\PXxG=\pi_1(\Xman,\Gxorb,\xelg)$ into a group with the following multiplication: $\pprod{[\agpath]}{[\bgpath]} := [\pconc{\agpath}{(\gpact{\ehomXx(\agpath)}{\bgpath})}]$.
Moreover, it directly follows from the definition that if $\ehomXx(\agpath) = \gelg$ and $\ehomXx(\bgpath) = \gelg'$, i.e.\ $\agpath(1)=\gact{\gelg}{\xelg}$ and $\bgpath(1)=\gact{\gelg'}{\xelg}$, then $(\pprod{\agpath}{\bgpath})(1) = \gact{(\gelg\gelg')}{\xelg}$.
In other words, $\ehomXx(\pprod{[\agpath]}{[\bgpath]})= \gprod{\ehomXx([\agpath])}{\ehomXx([\bgpath])}$, so $\ehomXx$ is a homomorphism, and its kernel is evidently $\pi_1(\Xman,\xelg)$.
This gives the short exact sequence~\eqref{equ:short_ex_seq_X_Gx_x}.

Statement~\ref{enum:lm:pi1Xax:iso_by_path} is also standard.

\ref{enum:lm:pi1Xax:equivariant_maps}
The left square of~\eqref{equ:func_morphism_of_sequences} is known to be commutative.
Suppose $\func$ is equivariant on $\Gxorb$.
We need to check commutativity of the right square of~\eqref{equ:func_morphism_of_sequences}.
Let $\agpath\in\PXxG$ and $\gelg =\ehomXx(\agpath)$, so $\agpath(1) = \gact{\ehomXx(\agpath)}{\xelg}$.
Then
\[(\func\circ\agpath)(1) = \func(\gact{\gelg}{\xelg})=\hact{\phi(\gelg)}{\yelg},\]
whence $\phi(\gelg) = \ehomYy(\func\circ\agpath)$.
In other words, $\phi\circ\ehomXx=\ehomYy\circ\func_{*}$.

Suppose that $\func$ is $\phi$-equivariant on all of $\Xman$, that is $\func\circ\gelg=\phi(\gelg)\circ\func\colon\Xman\to\Yman$ for all $\gelg\in\Ggrp$.
Let also $\agpath,\bgpath\in\PXxG$.
Then
\begin{align*}
    \func\circ\bigl(\pprod{\agpath}{\bgpath}\bigr)
    &= \func\circ\bigl(\pconc{\agpath}{(\gpact{\ehomXx(\agpath)}{\bgpath})}\bigr)
     = \pconc{(\func\circ\agpath)}{(\func\circ\gpact{\ehomXx(\agpath)}{\bgpath})} \\
    &= \pconc{(\func\circ\agpath)}{(\hpact{\phi(\ehomXx(\agpath))}{\func\circ\bgpath})}
     \stackrel{\eqref{equ:func_morphism_of_sequences}}{=\!=} \pconc{(\func\circ\agpath)}{(\hpact{\ehomYy(\func\circ\agpath)}{\func\circ\bgpath})} \\
    &= \pprod{(\func\circ\agpath)}{(\func\circ\bgpath)}.
\end{align*}
Hence, $\func_{*}(\pprod{[\agpath]}{[\bgpath]})=\pprod{[\func_{*}(\agpath)]}{[\func_{*}(\bgpath)]}$, so $\func_{*}$ is a homomorphism.
\end{proof}

Consider several particular cases of the constructions of Lemma~\ref{lm:pi1Xax}.
Suppose that we are given an action of $\Ggrp$ on $\Xman$ being WD at some $\xelg\in\Xman$.
We will regard elements of $\Ggrp$ as homeomorphisms of $\Xman$.

\begin{subexample}\rm
Let $\gelg\in\Ggrp$ and $\phi\colon \Ggrp\to\Ggrp$ be the inner automorphism of $\Ggrp$ induced by $\gelg$, i.e.\ $\phi(\kelg) = \gelg\circ\kelg\circ\gelg^{-1}$, $\kelg\in\Ggrp$.
Then the \term{homeomorphism} $\gelg\colon\Xman\to\Xman$ is $\phi$-equivariant:
\[
    \gelg\bigl( \kelg(\yelg) \bigr)
    = (\gelg \circ \kelg \circ \gelg^{-1})\bigl(\gelg(\yelg) \bigr)
    = \phi(\kelg)\bigl( \gelg(\yelg) \bigr), \qquad \kelg\in\Ggrp,\ \yelg\in\Xman.
\]
Hence, by Lemma~\ref{lm:pi1Xax}\ref{enum:lm:pi1Xax:equivariant_maps}, it induces an isomorphism of the following short exact sequences:
\begin{equation}\label{equ:conjugation}
    \begin{aligned}
    \xymatrix{
    \ \pi_1(\Xman,\xelg)                              \ \ar[d]^-{\gelg_{*}}  \ar@{^(->}[r] &
    \ \pi_1(\Xman,\Gxorb,\xelg)                       \ \ar[d]^-{\gelg_{*}}  \ar@{->>}[r]^-{\ehomXx} &
    \ \Ggrp                                           \ \ar[d]^-{\kelg\,\mapsto\gelg\circ\kelg\circ\gelg^{-1}} \\
    \ \pi_1(\Xman,\gact{\gelg}{\xelg})                \ \ar@{^(->}[r] &
    \ \pi_1(\Xman,\Gxorb,\gact{\gelg}{\xelg})         \ \ar@{->>}[r]^-{\ehomXgx} &
    \ \Ggrp                                           \
    }
\end{aligned}
\end{equation}
\end{subexample}

\begin{subexample}\rm
Let $\Kgrp$ be a normal subgroup of $\Ggrp$, so it also acts on $\Xman$, $\phi\colon \Ggrp\to\Ggrp/\Kgrp$ be the quotient homomorphism, $\Yman = \Xman/\Kgrp$ the quotient space endowed with the quotient topology, $\func\colon\Xman\to\Yman$ the quotient map, $\xelg\in\Xman$, and $\yelg = \func(\xelg)$.

Then $\Ggrp/\Kgrp$ naturally acts on $\Yman$.
Indeed, let $\yelg'\in\Yman$.
Then $\func^{-1}(\yelg') = \Kgrp\xelg'$ is the $\Kgrp$-orbit of some point $\xelg'\in\Xman$.
Now if $\gelg\in\Ggrp$, then
\begin{equation}\label{equ:gKx=Kgx}
    \gact{\gelg}{\Kgrp\xelg'}=\Kgrp\gact{\gelg}{\xelg'},
\end{equation}
since $\Kgrp$ is normal.
Hence, $\gact{\gelg}{\Kgrp\xelg'}$ depends only on the adjacent class $\phi(\gelg)$ of $\gelg$ in $\Ggrp/\Kgrp$, and thus we get an action of $\Ggrp/\Kgrp$ on $\Yman$.

Moreover, the identity~\eqref{equ:gKx=Kgx} can also be written as $\phi(\gelg)(\func(\xelg'))=\func(\gact{\gelg}{\xelg'})$ which means that $\func$ is $\phi$-equivariant.

Finally, the $\Ggrp/\Kgrp$-stabilizer of $\yelg$ is trivial, since so is the $\Ggrp$-stabilizer of $\xelg$.

Suppose, in addition, that the path component of $\yelg$ in its $\Ggrp/\Kgrp$-orbit consists of $\yelg$ only, that is the $\Ggrp/\Kgrp$-action on $\Yman$ is WD at $\yelg$.
Then, again by Lemma~\ref{lm:pi1Xax}\ref{enum:lm:pi1Xax:equivariant_maps}, we get a morphism of the following short exact sequences:
\begin{equation}\label{equ:qoutient}
    \begin{aligned}
    \xymatrix{
    \ \pi_1(\Xman,\xelg)                          \ \ar[d]^-{\func_{*}}  \ar@{^(->}[r] &
    \ \pi_1(\Xman,\Gxorb,\xelg)                   \ \ar[d]^-{\func_{*}} \ar@{->>}[r]^-{\ehomXx} &
    \ \Ggrp                                       \ \ar@{->>}[d]^-{\phi} \\
    \ \pi_1(\Xman/\Kgrp,\yelg)                    \ \ar@{^(->}[r] &
    \ \pi_1(\Xman/\Kgrp,\GKyorb,\yelg)            \ \ar@{->>}[r]^-{\ehomYy} &
    \ \Ggrp/\Kgrp                                 \
    }
\end{aligned}
\end{equation}
\end{subexample}

\subsection{Dynamical systems}\label{sect:dyn_syst}
Let $\Xman$ be a path connected topological space, $\xelg\in\Xman$ be a point, and $\phi\colon \Xman\to\Xman$ be a homeomorphism, so the iterations of $\phi$ define an action of $\bZ$ on $\Xman$.
If this action is WD at $\xelg$, then one can define the group $\pi_1(\Xman,\bZ\xelg,\xelg)$.
We will consider below few computations of that group, but first let us discuss the situations when such points exist.

\begin{subremark}\rm
a) Clearly, the condition that the stabilizer of $\xelg$ is trivial means that $\xelg$ is \term{non-periodic}, i.e.\ $\xelg\not=\phi^n(\xelg)$ for all $n\in\bZ\setminus 0$.

Note that existence of non-periodic points is a very typical situation.
Even more, for \qq{good} spaces like CW-complexes, homeomorphisms with only periodic points, and in particular, periodic homeomorphisms, are \qq{rare} in the corresponding homeomorphism groups: having all points periodic is an unstable property.
One might also mention a result by D.~Montgomery~\cite{Montgomery:AJM:1937} claiming that every homeomorphism of a \term{connected} manifold with all periodic points is itself periodic.

b) Assume further that $\xelg$ is non-periodic.
Then a $\bZ$-action is WD at $\xelg$ if $\{\xelg\}$ is the path components of $\xelg$ in its orbit $\bZ\xelg$.

This condition is also typical for \qq{good} spaces.
It may fail in some \qq{pathological} cases, e.g.\ when the orbit $\bZ\xelg$ or all $\Xman$ has anti-discrete topology (consisting only of two sets: $\Xman$ and $\varnothing$).

On the other hand, \term{if the one-point set $\{\xelg\}$ is closed in $\bZ\xelg$ (which holds e.g.\ when $\Xman$ is a $T_1$-space), then $\{\xelg\}$ is the path components of $\xelg$ in $\bZ\xelg$.}
Indeed, as noted above, since $\phi$ is a homeomorphism, for each $n\in\bZ$ the one-point set $\{\phi^n(\xelg)\}$ is also closed in $\bZ\xelg$.
Now let $\alpha\colon I\to\bZ\xelg$ be a path in the orbit, with $\alpha(0)=\xelg$.
Then for each $n\in\bZ$ the set $A_n = \alpha^{-1}(\phi^{n}(\xelg))$ is closed in $I$, and so we get at most countable partition $I = \sqcup_{n\in\bZ} A_n$ of $I$ into closed subsets.
By Sierpi\'{n}ski's theorem, (see e.g.~\cite{Cohen:APM:2013} for details), this is possible only when \term{all of these sets are empty except one of them}, which means that $A_0=I$, and so $\alpha$ is a constant path.
\end{subremark}

For instance, it follows from the above discussion that a $\bZ$-action is WD at each \term{wandering} point $\xelg$ of $\phi$.

\begin{subexample}\rm
Let $\alpha\in(0,1)$ and $\phi\colon \Circle\to\Circle$, $\phi(w)=we^{2\pi i \alpha}$, be the rotation of the circle.

1) Suppose $\alpha=a/m$ is rational, with $a,m\in\bN$ and $\mathrm{gcd}(a,m) = 1$.
Then $\phi$ generates a free action of $\bZ_m$ on $\Circle$.
One easily checks that for any $\welg\in\Circle$ the short exact sequence~\eqref{equ:short_ex_seq_X_Gx_x} for $(\Circle,\bZ_m\welg,\welg)$ is isomorphic with $m\bZ\monoArrow\bZ\epiArrow\bZ_m$.
In particular, $\pi_1(\Circle,\bZ_m\welg, \welg)=\bZ$.

2) If $\alpha$ is irrational, then $\phi$ generates a free action of $\bZ$.
One easily checks that now for any $\welg\in\Circle$, the short exact sequence~\eqref{equ:short_ex_seq_X_Gx_x} for $(\Circle,\bZ\welg, \welg)$ is isomorphic with
\[\bZ\oplus 0 \monoArrow \bZ^2 \epiArrow 0\oplus\bZ.\]
This will also follow from Lemma~\ref{lm:pi1XGxx} below, since $\phi$ is embeddable into a flow on $\Circle$, i.e.\ an action of $\bR$ and $\bZ$ is a subgroup of $\bR$.
\end{subexample}

\begin{subexample}\rm
Let $\phi\colon \Circle\to\Circle$, $\phi(\welg)=\overline{\welg}$, be the complex conjugation.
Then $\phi$ defines an action of $\bZ_2$ on $\Circle$ with two fixed points $\pm1$.
That action is WD at each $\welg\in\Circle\setminus\{\pm1\}$, and one checks that the short exact sequence~\eqref{equ:short_ex_seq_X_Gx_x} for $(\Circle,\bZ_2\welg, \welg)$ is $\bZ \monoArrow \bZ \tilde{\rtimes} \bZ_2 \epiArrow \bZ_2$ and corresponds to a unique non-trivial $\bZ_2$-extension $\bZ \tilde{\rtimes} \bZ_2$ of $\bZ$.
Recall that, by definition, $\bZ \tilde{\rtimes} \bZ_2$ is the semidirect product of $\bZ$ and $\bZ_2$ corresponding to the canonical isomorphism $\bZ_2 \xrightarrow{\cong} \Aut(\bZ)$.
More precisely, $\bZ \tilde{\rtimes} \bZ_2$ is the Cartesian product of sets $\bZ \times \bZ_2$ with the following operation: $(a,\delta) (b,\eps) := (a + \delta b, \delta\eps)$, where $a,b\in\bZ$ and $\delta,\eps\in\{\pm1\}=\bZ_2$.
One can also regard $\bZ \tilde{\rtimes} \bZ_2$ as the following group of integer $(2\times 2)$-matrices:
$ \bigl\{
    \bigl(
    \begin{smallmatrix}
       1 & 0 \\
       a & \delta
    \end{smallmatrix}
    \bigr)
    \mid
    a\in\bZ, \delta\in\{\pm1\}
 \bigr\}$.
\end{subexample}

Consider also one explicit computation of the sequence~\eqref{equ:short_ex_seq_X_Gx_x} for flows, i.e.\ actions of $\bR$.
Let $\flow\colon\Xman\times\bR\to\Xman$ be a flow on a topological space $\Xman$, and $\Ggrp\subset\bR$ be a subgroup distinct from $\bR$.
Let also $\xelg\in\Xman$ be a point and $\Gxorb = \{ \flow(\xelg,\tau) \mid \tau\in\Ggrp\}$ be its $\Ggrp$-orbit.

\begin{sublemma}\label{lm:pi1XGxx}\sl
Suppose that the path component of $\xelg$ in $\Gxorb$ is $\{\xelg\}$ and either
\begin{enumerate}[label={\rm(\alph*)}]
\item\label{enum:lm:pi1XGxx:nonper} $\xelg$ is non-periodic, or
\item\label{enum:lm:pi1XGxx:per}    $\xelg$ is periodic of some period $\theta>0$ such that $\theta\not\in\Ggrp$.
\end{enumerate}
Then the action of $\Ggrp$ is WD at $\xelg$, and the short exact sequence~\eqref{equ:short_ex_seq_X_Gx_x} splits, i.e.\ it is isomorphic with
\[
     \pi_1(\Xman,\xelg) \xmonoArrow{~\alpha\ \mapsto \ (\alpha, 0)}
     \pi_1(\Xman,\xelg) \times \Ggrp \xepiArrow{~(\alpha,\tau) \ \mapsto \ \tau~}
     \Ggrp.
\]
\end{sublemma}
\begin{proof}
Since $\Ggrp\subsetneq\bR$, the path component of $0\in\Ggrp$ must be $\{0\}$.
Then each of the conditions~\ref{enum:lm:pi1XGxx:nonper} and~\ref{enum:lm:pi1XGxx:per} imply that the stabilizer of $\xelg$ with respect to $\Ggrp$ is trivial.
As the path component of $\xelg$ in $\Gxorb$ is $\{\xelg\}$, we see that the action of $\Ggrp$ is WD at $\xelg$.

Let us compute the short exact sequence~\eqref{equ:short_ex_seq_X_Gx_x}:
$\pi_1(\Xman,\xelg) \xmonoArrow{j} \pi_1(\Xman,\Gxorb,\xelg) \xepiArrow{\ehomXx} \Ggrp$.
For each $\tau\in\Ggrp$ define the following path $\gamma_{\tau}\colon(I,0,1)\to(\Xman,\xelg,\Gxorb)$ by $\gamma_{\tau}(t) = \flow(x, t\tau)$, so it goes along the trajectory of $\xelg$ from $\xelg$ to $\flow(\xelg,\tau)$.

One easily check that $\ehomXx(\gamma_{\tau})=\tau$ and $\pconc{\gamma_{\tau}}{\gamma_{\tau'}}=\gamma_{\tau+\tau'}$ for all $\tau,\tau'\in\Ggrp$.
Therefore, $\Kgrp = \{[\gamma_{\tau}]\}_{\tau\in\Ggrp}$ is a subgroup of $\pi_1(\Xman,\Gxorb,\xelg)$ which is isomorphically mapped $\Kgrp$ onto $\Ggrp$ by $\ehomXx$.

It suffices to show that $\Kgrp$ commutes with $\pi_1(\Xman,\xelg)$.
This will imply that $\pi_1(\Xman,\Gxorb,\xelg)$ splits into the direct product of $\pi_1(\Xman,\xelg)$ and $\Kgrp$.

Let $\alpha\colon (I,\partial I) \to (\Xman,\xelg)$ be a loop at $\xelg$ and $\tau\in\Ggrp$.
Define the homotopy $F\colon I^2\to\Xman$ by $F(s,t)=\flow(\alpha(s), t\tau)$.
Denote by $\Jman = \overline{\partial I^2 \setminus (I\times 0)}$ the arc being the union of the left, top and right sides of $I^2$.
Then it is evident that $[\restr{F}{I\times 0}] = [\alpha]$, while $[\restr{F}{\Jman}] = [\pconc{\gamma_{\tau}}{\pconc{\alpha}\gamma_{\tau}^{-1}}]$.
Thus, $\alpha$ is homotopic to $\pconc{\gamma_{\tau}}{\pconc{\alpha}\gamma_{\tau}^{-1}}$ in $\pi_1(\Xman,\xelg)$, whence $\pi_1(\Xman,\xelg)$ commutes with $\Kgrp$.
\end{proof}

\section{Wreath products}\label{sect:wreath_products}
\subsection{Definitions}
Let $\ahm\colon\Xman\times\Kgrp\to\Xman$ be a \term{right} action of some group $\Kgrp$ on a set $\Xman$.
To simplify notation we will also write $\xelg\kelg$ instead of $\ahm(\xelg,\kelg)$.
For each $\kelg\in\Kgrp$ denote by $\RightShift{\kelg}\colon\Xman\to\Xman$, $\RightShift{\kelg}(\xelg)=\xelg\kelg$, the \term{shift} of $\Xman$ by $\kelg$.
Evidently, $\RightShift{\lelg}\circ\RightShift{\kelg}(\xelg)=\RightShift{\lelg}(\xelg\kelg)=\xelg\kelg\lelg=\RightShift{\kelg\lelg}(\xelg)$.

Let also $\Ggrp$ be another group.
Then the set $\Maps{\Xman}{\Ggrp}$ of all maps $\Xman\to\Ggrp$ is a group with respect to the point-wise multiplication, and we also have a natural \term{left} action
\[
    \kappa\colon\Kgrp\times\Maps{\Xman}{\Ggrp}\to \Maps{\Xman}{\Ggrp},
    \qquad
    \kappa(\kelg,\amapgf) = \amapgf \circ \RightShift{\kelg}\colon \Xman \xrightarrow{\RightShift{\kelg}} \Xman \xrightarrow{\alpha} \Ggrp,
\]
of $\Kgrp$ on $\Maps{\Xman}{\Ggrp}$.
It is indeed a \term{left} action, since
\[
    \kappa( \kelg, \kappa(\lelg,\amapgf) ) =
    \kappa( \kelg, \amapgf \mcirc \RightShift{\lelg}) =
    \amapgf \mcirc \RightShift{\lelg}\mcirc\RightShift{\kelg} =
    \amapgf \mcirc \RightShift{\kelg\lelg}.
\]

The semidirect product $\Maps{\Xman}{\Ggrp} \rtimes_{\ahm} \Kgrp$ associated with this action is called the \emph{(unrestricted) wreath product} of $\Ggrp$ and $\Kgrp$ corresponding to $\ahm$ and will be denoted by $\Ggrp\wrm{\ahm}\Kgrp$ or even $\Ggrp\wrm{\ahm,\Xman}\Kgrp$ if we need to specify $\Xman$.
More precisely, the multiplication in $\Ggrp\wrm{\ahm,\Xman}\Kgrp$ is defined as follows: if $\alpha,\beta\colon\Xman\to\Ggrp$ are two maps and $\kelg,\lelg\in\Kgrp$, then
\[
    (\alpha,\kelg)(\beta,\lelg) := (\gprod{\alpha}{(\beta\mcirc\RightShift{\kelg})},\kelg\lelg),
\]
where $\cdot$ means the point-wise multiplication of maps, i.e.\ $(\gprod{\alpha}{(\beta\mcirc\RightShift{\kelg})})(\xelg) = \alpha(\xelg)\beta(\xelg\kelg)$ for all $\xelg\in\Xman$.

Let $\unit{\Ggrp}$ and $\unit{\Kgrp}$ be the units of $\Ggrp$ and $\Kgrp$ respectively, and $\mathbf{e}_{\Ggrp}\colon\Xman\to\Ggrp$ be the constant map into the unit $\unit{\Ggrp}\in\Ggrp$.
Then $(\mathbf{e}_{\Ggrp}, \unit{\Kgrp})$ is the unit of $\Ggrp\wrm{\ahm,\Xman}\Kgrp$, and the inverse of $(\alpha,\kelg)$ is $((\alpha\mcirc\RightShift{\kelg^{-1}})^{-1}, \kelg^{-1})$, where $(\alpha\mcirc\RightShift{\kelg^{-1}})^{-1}\colon\Xman\to\Ggrp$ is the point-wise inverse of $\alpha\mcirc\RightShift{\kelg^{-1}}\colon\Xman\to\Ggrp$, so $(\alpha\mcirc\RightShift{\kelg^{-1}})^{-1}(\xelg) = (\alpha(\xelg\kelg^{-1}))^{-1}$ for all $\xelg\in\Xman$.

Again, the \term{restricted wreath product} is defined similarly, but one should replace $\Maps{\Xman}{\Ggrp}$ with its subset $\MapsFin{\Xman}{\Ggrp}$ of functions $\alpha\colon \Xman\to\Ggrp$ whose \term{support}, $\Xman\setminus\alpha^{-1}(\unit{\Ggrp})$, is finite.

We will be interested in the following two particular cases of this construction for the situation when $\Xman$ itself is a group and $\Kgrp$ acts on $\Xman$ by left shifts.

\subsection{Regular wreath products of groups corresponding to effective actions}
Suppose $\Xman = \Kgrp$ and the action $\ahm\colon\Kgrp\times\Kgrp\to\Kgrp$ is just the multiplication in $\Kgrp$.
In this case
\[ \Ggrp\wrm{\ahm,\Kgrp}\Kgrp := \Maps{\Kgrp}{\Ggrp}\rtimes_{\ahm}\Kgrp \]
is denoted simply by $\Ggrp\wr\Kgrp$ and usually called the \term{regular wreath product} of $\Ggrp$ and $\Kgrp$.

\begin{subexample}\rm
Let $\Kgrp$ be a finite group of some order $m\geq1$.
Then $\Maps{\Kgrp}{\Ggrp}$ can be identified with the $m$-th Cartesian power $\Ggrp^{m}$ of $\Ggrp$, whose coordinates are enumerated by elements of $\Kgrp$, and
$\Ggrp\wr\Kgrp$ is the product of sets $\Ggrp^{m}\times\Kgrp$ with the following multiplication:
\begin{equation}\label{equ:GwrK_mult}
\bigl( \{a_{\kelg}\}_{\kelg\in\Kgrp}, \gelg \bigr)
\bigl( \{b_{\kelg}\}_{\kelg\in\Kgrp}, \helg \bigr) =
\bigl( \{a_{\kelg}b_{\kelg\gelg}\}_{\kelg\in\Kgrp}, \gelg\helg \bigr)
\end{equation}
where $a_i,b_j\in\Ggrp$, $k,l\in\Kgrp$.
\end{subexample}

\begin{subexample}\rm
In particular, if $\Kgrp=\bZ_{m}$ is a finite cyclic group of some order $m\geq1$, then $\Ggrp\wr\bZ_{m}$ can be regarded as the product of \term{sets} $\Ggrp^{m}\times\bZ_m$ with the following multiplication:
\begin{equation}\label{equ:GwrZm_mult}
(a_0,\ldots,a_{m-1};\pelg) (b_0,\ldots,b_{m-1};\qelg) =
(a_{0} b_{\pelg}, a_{1} b_{1+\pelg}, \ldots, a_{m-1} b_{\pelg-1}, \pelg + \qelg),
\end{equation}
where $a_i,b_j\in\Ggrp$, $\pelg,\qelg\in\bZ_{m}$, and all indices are taken modulo $m$.
\end{subexample}

\begin{subexample}\rm
\label{enum:W4}
Similarly, let $\Kgrp=\bZ_{m}\times\bZ_{n}$ be the product of two finite cyclic groups of orders $m,n\geq1$.
Then $\Maps{\bZ_{m}\times\bZ_{n}}{\Ggrp}$ can be identified with $mn$-th Cartesian power $\Ggrp^{mn}$ of $\Ggrp$, whose elements can be regarded as $(m\times n)$-matrices with entries in $\Ggrp$.
Moreover, $\Ggrp\wr(\bZ_{m}\times\bZ_{n})$ is the product of sets $\Ggrp^{mn}\times\bZ_{m}\times\bZ_{n}$ with the following multiplication:
\begin{multline}\label{equ:GwrZmZn_mult}
\bigl( \{a_{i,\,j}\}_{i\in\bZ_m,\,j\in\bZ_n}; p, q \bigr) \,
\bigl( \{b_{i,\,j}\}_{i\in\bZ_m,\,j\in\bZ_n}; r, s \bigr) \, = \\ = \,
\bigl( \{a_{i,\,j}\, b_{i+p,\,j+q} \}_{i\in\bZ_m,\,j\in\bZ_n}; p+r, q+s) \bigr),
\end{multline}
where $a_{i,\,j}, b_{i,\,j}\in\Ggrp$, $p,r\in\bZ_{m}$, $q,s\in\bZ_{n}$, and all first and second indices are taken modulo $m$ and $n$ respectively.
\end{subexample}

\subsection{Wreath products corresponding to non-effective actions}
\label{exmp:wreath_prod_of_homomorphism}
There is an extensive literature on wreath products, see e.g.~\cite{Wells:AMM:1976, Meldrum:Wreath:1995, BhattacharjeeMacphersonMollerNeumann:LNM:1997} and references therein.
Most of them correspond to effective actions, see also~\cite{Hold:QJMS:1978}.
However, it was recently proved by the author, \cite{Maksymenko:TA:2020}, that for typical smooth functions $\func$ on with isolated critical points of compact surfaces $\Mman$, the fundamental groups of their orbits with respect to natural actions of the diffeomorphism groups $\Diff(\Mman)$ are ``built'' from wreath products corresponding to certain non-effective actions of $\bZ$, see Sections~\ref{sect:classes_B},\ref{sect:orbits}.
We will define below those ``building blocks''.

Assume now that $\Xman$ is a group and let $\hahm\colon\Kgrp\to\Xman$ be a homomorphism.
Then we have a natural right action $\ahm\colon\Xman\times\Kgrp\to\Xman$ of $\Kgrp$ on $\Xman$ (as a set) by right shifts given by $\ahm(\xelg,\kelg) = \xelg\hahm(\kelg)$.
In this case the corresponding wreath product $\Ggrp\wrm{\ahm,\Xman}\Kgrp$ will also be denoted by $\Ggrp\wrm{\hahm\colon\Kgrp\to\Xman}\Kgrp$ or simply by $\Ggrp\wrm{\hahm}\Kgrp$.
Again, it is the set $\Maps{\Xman}{\Ggrp}\times\Kgrp$ with the following multiplication: if $\alpha,\beta\colon\Xman\to\Ggrp$ are two maps and $\kelg,\lelg\in\Kgrp$, then
\[
    (\alpha,\kelg)(\beta,\lelg) := (\gprod{\alpha}{\beta\mcirc\RightShift{\hahm(\kelg)}},\kelg\lelg).
\]
Notice that if $\hahm$ has a non-trivial kernel, then the corresponding action $\ahm$ is \term{non-effective}.

Evidently, the regular wreath product $\Ggrp\wr\Kgrp$ is the same as $\Ggrp\wrm{\id_{\Kgrp}}\Kgrp$, i.e.\ it corresponds to the identity isomorphism of $\Kgrp$.

\begin{sublemma}\sl
Let $\Kgrp,\Xman$ and $\Yman$ be any groups.
Then every homomorphism $\phi\colon \Xman\to\Yman$ induces a homomorphism
\[
\phi^{*}\colon\
    \Ggrp\wrm{\phi\circ\hahm\colon\Kgrp\to\Yman}\Kgrp
    \ \longrightarrow \
    \Ggrp\wrm{\hahm\colon\Kgrp\to\Xman}\Kgrp,
    \qquad
    \phi^{*}(\alpha,\kelg) = (\alpha\circ\phi,\kelg),
\]
for $\alpha\colon \Yman\to\Ggrp$ and $\kelg\in\Kgrp$.
Moreover, $(\id_{\Xman})^{*} = \id_{(\Ggrp\wrm{\hahm\colon\Kgrp\to\Xman}\Kgrp)}$, and if $\psi\colon\Yman\to\Zman$ is another homomorphism of groups, then $(\psi\circ\phi)^{*} = \phi^{*}\circ\psi^{*}$.
In particular, if $\phi$ is an isomorphism, then so is $\phi^{*}$.
\end{sublemma}
\begin{proof}
Note that for every $\kelg,\kelg'\in\Kgrp$ we have that
\begin{equation}\label{equ:R__phi_commute}
    \RightShift{\phi(\hahm(\kelg))}\circ \phi(\kelg') =
    \phi(\kelg') \, \phi(\hahm(\kelg)) =
    \phi(\kelg' \hahm(\kelg)) =
    \phi\circ \RightShift{\hahm(\kelg)}(\kelg').
\end{equation}
Hence, if $(\alpha,\kelg), (\beta,\lelg) \in \Ggrp\wrm{\phi\circ\hahm\colon\Kgrp\to\Yman}\Kgrp = \Maps{\Yman}{\Ggrp}\rtimes\Kgrp$, then
\begin{align*}
\phi^{*}\bigl( (\alpha,\kelg)(\beta,\lelg) \bigr) &=
\bigl(
    (\gprod{\alpha}{\beta\mcirc\RightShift{\phi(\hahm(\kelg))}}) \circ \phi, \kelg\lelg
\bigr)
=
\bigl(
    \gprod{(\alpha\circ \phi)}{(\beta\mcirc\RightShift{\phi(\hahm(\kelg))}\circ \phi)}, \kelg\lelg
\bigr) \\
&\stackrel{\eqref{equ:R__phi_commute}}{=}
\bigl(
    \gprod{(\alpha\circ \phi)}{(\beta\mcirc\phi\circ\RightShift{\hahm(\kelg)})}, \kelg\lelg
\bigr) =
\phi^{*}(\alpha,\kelg) \phi^{*}(\beta,\lelg).
\end{align*}
Thus, $\phi^{*}$ is a homomorphism.
All other statements are evident.
\end{proof}

\begin{subexample}\rm
\label{enum:W6}
Let $\Kgrp=\bZ$, $\Xman=\bZ_{m}$ and $\hahm\colon\bZ\to\bZ_{m}$, $\hahm(k)=k\bmod m$, be the natural $\bmod\,m$ epimorphism.
Then the group $\Ggrp\wrm{\hahm\colon\bZ\to\bZ_{m}} \bZ$ will be denoted by $\Ggrp\wrm{m}\bZ$.
It can be regarded as the product of sets $\Ggrp^{m}\times\bZ$ with the multiplication given precisely by the same formula~\eqref{equ:GwrZm_mult} in which $\pelg,\qelg$ now belong to $\bZ$.
\end{subexample}

\begin{subexample}\rm
\label{enum:W7}
Similarly, let $\hahm\colon\bZ^2\to\bZ_{m}\times\bZ_{n}$, $\hahm(k,l)=(k\bmod m, l\bmod n)$.
Then the group $\Ggrp\wrm{\hahm} \bZ$ will be denoted by $\Ggrp\wrm{m,n}\bZ^2$.
It can be regarded as the product of sets $\Ggrp^{mn}\times\bZ^2$ with the multiplication given precisely by the same formula~\eqref{equ:GwrZmZn_mult} in which $p,q,r,s$ now belong to $\bZ$.
\end{subexample}

\begin{subexample}\rm
\label{enum:W8}
More generally, let $\Kgrp=\bZ_{m_1}\times\cdots\times\bZ_{m_n}$ be the product of $n$ finite cyclic groups of orders $m_1,\ldots,m_n\geq1$, and $m=m_1 m_2 \cdots m_n$ be their product.
Define also the homomorphism $\hahm\colon\bZ^n\to\Kgrp$, $\hahm(k_1,\ldots,k_n)=(k_1\bmod m_1, \ldots, k_n\bmod m_n)$.
Again for any group $\Ggrp$ the set $\Maps{\Kgrp}{\Ggrp}$ can be identified with $m$-th Cartesian power $\Ggrp^{m}=G^{m_1}\times\cdots \times G^{m_n}$ of $\Ggrp$.
Then similarly to~\ref{enum:W4} and~\ref{enum:W7} one can define two wreath products: $\Ggrp\wr\Kgrp$ and $\Ggrp\wrm{\ahm}\bZ^{n}$.
They correspond respectively to the effective action of $\Kgrp$ and the non-effective action of $\bZ^{n}$ on $\Ggrp^{m_1}\times\cdots\times\Ggrp^{m_n}$ by independent cyclic shifts of coordinates.
We will also denote \[ \Ggrp\wrm{\ahm}\bZ^{n}:=\Ggrp\wrm{m_1,\ldots,m_n}\bZ^{n}.\]
\end{subexample}

\subsection{Short exact sequences}\label{sect:short_ex_sequence}
We will need to consider several operations over short exact sequences related with previous examples.
First define the following short exact sequences:
\begin{align}\label{equ:ex_seq_TZ}
\seqTriv& \colon 1\monoArrow 1 \epiArrow 1, &
\seqC{1}& \colon \bZ\xmonoArrow{\id} \bZ \epiArrow 1, &
\seqC{m}& \colon m\bZ\xmonoArrow{~~} \bZ \xepiArrow{~\mathrm{mod}\, m~} \bZ_m,
\end{align}
for $m\geq2$.
Further, given two short exact sequences $\pSeq\colon A\monoArrow B \epiArrow C$ and $\qSeq\colon K\monoArrow L \epiArrow M$, one can define their product
\[
    \pSeq\times\qSeq\colon A\times K\monoArrow B\times L \epiArrow C\times M.
\]
Also, for a short exact sequence $\aSeq\colon \kA\xmonoArrow{\alpha} \kB \xepiArrow{~\beta~} \kC$ and $m\geq1$ we have the following exact $(3\times3)$-diagram:
\begin{equation}\label{equ:wr_ex_seq_m}
\aligned
\xymatrix@C=1.8em@R=1.8em{
\ \aSeq^{m}\colon \ar@{^(->}[d] &
\kA^m\times0 \ar@{^(->}[d] \ar@{^(->}[r]  &
\kB^m\times0 \ar@{^(->}[d]  \ar@{->>}[r]  &
\kC^m\times0 \ar@{^(->}[d] \\
\ \seqWrm{\aSeq}{m}\colon \ar@{->>}[d] &
\kA^m\times m\bZ \ \ar@{^(->}[r]^-{\alpha'} \ar@{->>}[d] &
\kB\wrm{m}\bZ \ar@{->>}[r]^-{\beta'} \ar@{->>}[d]^-{\,p} &
\kC\wr\bZ_m \ar@{->>}[d] \\
\ \seqC{m}\colon &
m\bZ \ar@{^(->}[r] &
\bZ \ar@{->>}[r]  &
\bZ_m
}
\endaligned
\end{equation}
where $p\colon\kB\wrm{m}\bZ = \kB^m\rtimes\bZ\to\bZ$ is the projection onto the last coordinate, and
\begin{gather*}
\alpha'(a_1,\ldots,a_m, mk) = (\alpha(a_1),\ldots,\alpha(a_m), mk), \\
\beta'(b_1,\ldots,b_m, l) = (\beta(b_1),\ldots,\beta(b_m), \BMOD{l}{m}),
\end{gather*}
$a_i\in \kA$, $b_i\in \kB$ for all $i=1,\ldots,m$, and $k,l\in\bZ$.
The middle horizontal sequence will be denoted by $\seqWrm{\aSeq}{m}$.
Thus, \eqref{equ:wr_ex_seq_m} can be viewed as a short exact sequence of its rows: $\aSeq^{m}\monoArrow\seqWrm{\aSeq}{m}\epiArrow\seqC{m}$.

More generally, let $m_1,\ldots,m_k\geq1$ be natural numbers, and $m=m_1\cdots m_k$.
Then one has the following exact $(3\times 3)$-diagram:
\begin{equation}\label{equ:wr_ex_seq_m1mk}
\aligned
\xymatrix@C=1.8em@R=1.8em{
\ \pSeq^{m}\colon \ar@{^(->}[d] &
\kA^m\times0 \ar@{^(->}[d] \ar@{^(->}[r]  &
\kB^m\times0 \ar@{^(->}[d]  \ar@{->>}[r]  &
\kC^m\times0 \ar@{^(->}[d] \\
\ \seqWrm{\pSeq}{m_1,\ldots,m_k}\colon \ar@{->>}[d] &
\kA^m\times \prod\limits_{i=1}^{k} m_i\bZ \ \ar@{^(->}[r]^-{\alpha'} \ar@{->>}[d] &
\kB\wrm{m_1,\ldots,m_k}\bZ^k \ar@{->>}[r]^-{\beta'} \ar@{->>}[d]^-{\,p} &
\kC\wr\prod\limits_{i=1}^{k} \bZ_{m_i} \ar@{->>}[d] \\
\ \prod\limits_{i=1}^{k}\seqC{m_i}\colon &
\prod\limits_{i=1}^{k} m_i\bZ \ar@{^(->}[r] &
\bZ^k \ar@{->>}[r]  &
\prod\limits_{i=1}^{k} \bZ_{m_i}
}
\endaligned
\end{equation}
where $\alpha'$ and $\beta'$ defined in a similar way, and the middle horizontal sequence denoted by $\seqWrm{\pSeq}{m_1,\ldots,m_k}$.
Again, \eqref{equ:wr_ex_seq_m1mk} can be regarded as a short exact sequence
\[
\pSeq^{m}\monoArrow\seqWrm{\pSeq}{m_1,\ldots,m_k}\epiArrow\prod\limits_{i=1}^{k}\seqC{m_i}
\]
of its rows.

\section{Actions of wreath products}\label{sect:act_wreath_prod}
Suppose that we are given a \term{left} action of a group $\Ggrp$ on a set $\Xman$ and a \term{left} action of a group $\Hgrp$ on a set $\Yman$.
Let also $\Maps{\Hgrp}{\Xman}$ be the set of all maps $\Hgrp\to\Xman$ and $\Wman = \Maps{\Hgrp}{\Xman} \times \Yman$.

For each $\helg\in\Hgrp$ denote by $\RightShift{\helg}\colon\Hgrp\to\Hgrp$, $\RightShift{\helg}(\kelg)=\kelg\helg$, the \term{right} shift of $\Hgrp$ by element $\helg$.
The following statement is well-known.
We recall precise formulas in the proof, since they will be used in our main result.

\begin{sublemma}[{cf.~\cite[Sect.~2.3.1]{Arnold:PhD:2013}}]\label{lm:wreath_product_action}\sl
\begin{enumerate}[wide, itemsep=1ex, topsep=1ex, label={\rm\arabic*)}]
\item\label{enum:wr_act:def}
The wreath product $\Ggrp \wr \Hgrp$ acts from the left on $\Wman$ by the following rule: if $(\agfunc\colon\Hgrp\to\Ggrp, \helg) \in \GwrH = \Maps{\Hgrp}{\Ggrp}\rtimes\Hgrp$ and $(\xfunc\colon\Hgrp\to\Xman,\yelg) \in \Wman$, then
\begin{equation}\label{equ:wreath_prod_action}
    \ghact{(\agfunc,\helg)}{(\xfunc, \yelg)} =
    \bigl(\gprod{\agfunc}{\xfunc\mcirc\RightShift{\helg}}, \hact{\helg}{\yelg}\bigr),
\end{equation}
where $\gprod{\agfunc}{\xfunc\mcirc\RightShift{\helg}}\colon \Hgrp\to\Xman$ is the point-wise (in $\Hgrp$) multiplication, i.e.\ the map given by
\[
    (\gprod{\agfunc}{\xfunc\mcirc\RightShift{\helg}})(\kelg) = \gprod{\agfunc(\kelg)}{\xfunc(\kelg\helg)},
    \qquad
    \kelg\in\Hgrp.
\]

\item\label{enum:wr_act:stab}
Let $\xelg\in\Xman$, $\yelg\in\Yman$, $\cpathx\colon\Hgrp\to \{\xelg\}\subset \Xman$ be a constant map into the point $\xelg$, and $\welg = (\cpathx, \yelg) \in \Wman$.
Suppose that the $\Ggrp$-stabilizer of $\xelg$ and the $\Hgrp$-stabilizer of $\yelg$ are trivial.
Then the $(\GwrH)$-stabilizer of $\welg$ is also trivial, and
\begin{equation}\label{equ:GHobh_Gorb_Horb}
    \GHworb = \Maps{\Hgrp}{\Gxorb}\times\Hyorb.
\end{equation}

\item\label{enum:wr_act:free}
If the actions of $\Ggrp$ and $\Hgrp$ are free, then the action of $\GwrH$ is also free.
\end{enumerate}
\end{sublemma}
\begin{proof}
\ref{enum:wr_act:def}
First let us check that~\eqref{equ:wreath_prod_action} is a \term{left} action.
Let $(\xfunc,\yelg) \in \Wman$, and $(\agfunc, \kelg), (\bgfunc, \lelg) \in \GwrH$.
Recall that $(\agfunc, \kelg) (\bgfunc, \lelg) = (\gprod{\agfunc}{\bgfunc\mcirc\RightShift{\kelg}}, \kelg\lelg)$.
Hence,
\begin{align*}
\ghact{(\agfunc, \kelg)}{\bigl( \ghact{(\bgfunc,\lelg)}{(\xfunc, \yelg)} \bigr) }
&= \ghact{(\agfunc, \kelg)}
         { \bigl(
                \gprod{\bgfunc}{\xfunc\mcirc\RightShift{\lelg}},
                \hact{\lelg}{\yelg})
            \bigr)
         } \\
&= \bigl(
    \gprod{\agfunc}{(\gprod{\bgfunc}{\xfunc\mcirc\RightShift{\lelg}})\mcirc\RightShift{\kelg}},
    \hact{\kelg}{(\hact{\lelg}{\yelg})}
    \bigr) \\
&= \bigl(
    \gprod{\agfunc}{\gprod{\bgfunc\mcirc\RightShift{\kelg}}{\xfunc\mcirc\RightShift{\lelg}\mcirc\RightShift{\kelg}}},
    \hact{(\kelg\lelg)}{\yelg}
    \bigr) \\
&= \bigl(
    \gprod{\agfunc}{\gprod{\bgfunc\mcirc\RightShift{\kelg}}{\xfunc\mcirc\RightShift{\kelg\lelg}}},
    \hact{(\kelg\lelg)}{\yelg}
    \bigr) =
    \ghact{\bigl( (\agfunc, \kelg)(\bgfunc,\lelg) \bigr)}{(\xfunc, \yelg)}.
\end{align*}

\ref{enum:wr_act:stab}
Suppose that the $\Ggrp$-stabilizer of $\xelg$ and the $\Hgrp$-stabilizer of $\yelg$ are trivial. Let also $(\agfunc,\helg)\in\GwrH$ be an element of the $\GwrH$-stabilizer of $\welg=(\cpathx, \yelg)$.
We should show that $(\agfunc,\helg) = (\mathbf{e}_{\Ggrp}, \unit{\Hgrp})$ is the unit of $\GwrH$.
Indeed, since $\cpathx$ is a constant map, $\cpathx = \cpathx\mcirc\RightShift{\helg}$, whence the relation $\welg = \ghact{(\agfunc,\helg)}{\welg}$ can be written as follows:
\begin{equation}\label{equ:ah__on__xy}
(\cpathx, \yelg) = \ghact{(\agfunc,\helg)}{(\cpathx, \yelg)}
    = \bigl(\gprod{\agfunc}{\cpathx\mcirc\RightShift{\helg}}, \hact{\helg}{\yelg}\bigr)
    = \bigl(\gprod{\agfunc}{\cpathx}, \hact{\helg}{\yelg}\bigr).
\end{equation}
Comparing first coordinates of~\eqref{equ:ah__on__xy}, we see that $\cpathx = \gprod{\agfunc}{\cpathx}$ which means that for every $\kelg\in\Hgrp$,
\[
    \xelg = \cpathx(\xelg)
          = \gprod{\agfunc}{\cpathx}(\kelg)
          = \gact{\agfunc(\kelg)}{\xelg}.
\]
As the stabilizer of $\xelg$ is trivial, it then follows that $\agfunc(\kelg)=\unit{\Ggrp}$, so $\agfunc=\mathbf{e}_{\Ggrp}$.

Similarly, from the equality of second coordinates in~\eqref{equ:ah__on__xy} we obtain that $\yelg=\hact{\helg}{\yelg}$.
As the $\Hgrp$-stabilizer of $\yelg$ is trivial, we must have that $\helg=\unit{\Hgrp}$.
Thus, $(\agfunc,\helg) = (\mathbf{e}_{\Ggrp}, \unit{\Hgrp})$.

\medskip
It remains to prove that $\GHworb = \Maps{\Hgrp}{\Gxorb}\times\Hyorb$.
The identity~\eqref{equ:ah__on__xy} shows that $\GHworb \subset \Maps{\Hgrp}{\Gxorb}\times\Hyorb$.
Conversely, let $(\bgfunc,\yelg')\in\Maps{\Hgrp}{\Gxorb}\times\Hyorb$.
Since the $\Ggrp$-stabilizer of $\xelg$ is trivial, we see that for each $\kelg\in\Hgrp$ there exists a unique $\agfunc(\kelg)\in\Ggrp$ such that $\bgfunc(\kelg)=\gact{\agfunc(\kelg)}{\xelg}$, and thus $\bgfunc = \gprod{\agfunc}{\cpathx}$.
Also, since the $\Hgrp$-stabilizer of $\yelg$ is also trivial, there exists a unique $\helg\in\Hgrp$ such that $\yelg' = \hact{\helg}{\yelg}$.
Thus, $(\bgfunc,\yelg') = \bigl( \gprod{\agfunc}{\cpathx}, \hact{\helg}{\yelg} \bigr)=\ghact{(\agfunc,\helg)}{\welg}\in \GHworb$.

\ref{enum:wr_act:free}
Suppose the actions of $\Ggrp$ and $\Hgrp$ are free.
Assume also that
\[
    (\xfunc, \yelg)
    = \ghact{(\agfunc,\helg)}{(\xfunc, \yelg)}
    = (\gprod{\agfunc}{\xfunc\mcirc\RightShift{\helg}}, \hact{\helg}{\yelg})
\]
for some $(\agfunc,\helg) \in \GwrH$ and $(\xfunc, \yelg) \in \Wman$.
In particular, $\yelg=\hact{\helg}{\yelg}$, and since the action of $\Hgrp$ is free, we obtain that $\helg=\unit{\Hgrp}$.
Moreover, we also have that $\xfunc = \gprod{\agfunc}{\xfunc\mcirc\RightShift{\helg}} = \gprod{\agfunc}{\xfunc}$.
That is $\xfunc(\kelg)=\gprod{\agfunc(\kelg)}{\xfunc(\kelg)}$ for all $\kelg\in\Hgrp$.
Since the action of $\Ggrp$ on $\Xman$ is also free, it follows that $\agfunc(\kelg)=\unit{\Ggrp}$ for all $\kelg\in\Hgrp$.
In other words, $\agfunc = \mathbf{e}_{\Ggrp}\colon\Hgrp\to\Ggrp$ is the constant map into the unit of $\Ggrp$, whence $(\agfunc,\helg) = (\mathbf{e}_{\Ggrp}, \unit{\Hgrp})$ is the unit of $\GwrH$.
\end{proof}

\subsection{WD-actions of wreath products}
Assume now that $\Xman$ and $\Yman$ are path connected topological spaces, and the groups $\Ggrp$ and $\Hgrp$ act on them by homeomorphisms.
Fix two points $\xelg\in\Xman$ and $\yelg\in\Yman$.
Let also $\cpathx\colon\Hgrp\to \{\xelg\}\subset \Xman$ be a constant map into the point $\xelg$ and
\begin{equation}\label{equ:point_w}
    \welg = (\cpathx, \yelg) \in \Wman.
\end{equation}

Suppose also that $\Hgrp$ is finite of some order $m$.
Endow $\Hgrp$ with the discrete topology.
Then every map $\Hgrp\to\Xman$ is continuous, i.e.\ $\Maps{\Hgrp}{\Xman} = \Cont{\Hgrp}{\Xman}$, and we endow this space with the compact open topology.

Fix some enumeration $\Hgrp = \{\helg_1, \ldots, \helg_m \}$ of elements of $\Hgrp$.
Then the natural identification
\begin{align*}
   &\omega: \Maps{\Hgrp}{\Xman} \to \Xman^{m},
   &
   &\omega(\func) = (\func(\helg_1), \ldots,\func(\helg_m)),
\end{align*}
is a homeomorphism.
In particular, we get a homeomorphism
\[
    \omega\times\id_{\Yman}\colon \Wman =\Maps{\Hgrp}{\Xman} \times \Yman \xrightarrow{\omega\times\id_{\Yman}} \Xman^{m} \times \Yman.
\]
In what follows it will be convenient to regard $\Wman$ sometimes as $\Maps{\Hgrp}{\Xman}$ and sometimes as $\Xman^{m} \times \Yman$.
This will simplify some formulas and will never lead to confusion.

\begin{sublemma}\label{lm:GwrH_W_PD}\sl
\begin{enumerate}[leftmargin=*, label={\rm(\arabic*)}, itemsep=0.5ex, topsep=0.5ex]
\item\label{enum:GwrH:WD}
Suppose the action of $\Ggrp$ is WD at $\xelg$, while the action of $\Hgrp$ is WD at $\yelg$, (see Definition~\ref{def:WD}).
Then the action of $\GwrH$ on $\Wman$ is WD at $\welg$.

\item\label{enum:GwrH:PD}
If the actions of $\Ggrp$ and $\Hgrp$ are PD, then the action of $\GwrH$ on $\Wman$ is PD as well.
\end{enumerate}
\end{sublemma}
\begin{proof}
We regard $\Wman$ here as $\Xman^{m}\times\Yman$.

\ref{enum:GwrH:WD}
By assumption the $\Ggrp$-stabilizer of $\xelg$ and the $\Hgrp$-stabilizer of $\yelg$ are trivial, and every continuous path into $\Gxorb$ as well as into $\Hyorb$ is constant.
Then by Lemma~\ref{lm:wreath_product_action} the $\GwrH$-stabilizer of $\welg$ is also trivial.
Moreover, by~\eqref{equ:GHobh_Gorb_Horb}, $\GHworb = (\Gxorb)^{m}\times\Hyorb$, whence every continuous path into $\GHworb$ is also constant.
Hence, the action of $\GwrH$ is WD at $\welg$.

\ref{enum:GwrH:PD}
Let $(\xelg_1,\ldots,\xelg_m,\yelg_0)\in\Wman$ be any point.
Since the action of $\Hgrp$ on $\Yman$ is PD, there exists an $\Hgrp$-wandering neighborhood of $\Vman$ of $\yelg_0$, i.e.\ $\hact{\kelg}{\Vman} \cap \hact{\lelg}{\Vman} = \varnothing$ for all $\kelg\not=\lelg\in\Hgrp$.
Moreover, as the $\Ggrp$-action on $\Xman$ is also PD, for each $\xelg_i\in\Xman$, $i=1,\ldots,m$, there exists a $\Ggrp$-wandering neighborhood $\Uman_{i}$.
We claim that then $\Wman := \Uman_1\times\cdots\times\Uman_m\times\Vman$ is a $\GwrH$-wandering neighborhood of $(\xelg_1,\ldots,\xelg_m,\yelg_0)$ in $\Wman$.

Indeed, suppose there exists $\sigma = (\agfunc_1,\ldots,\agfunc_m,\helg)\in\GwrH$ such that $\Wman \cap \sigma(\Wman) \not=\varnothing$.
In other words, there exists $(\xelg'_1,\ldots,\xelg'_m,\yelg')\in \Wman$ such that
\[
    \ghact{(\agfunc_1,\ldots,\agfunc_m,\helg)}{(\xelg'_1,\ldots,\xelg'_m,\yelg')} =
    (\gact{\agfunc_{j_1}}{\xelg'_1}, \ldots, \gact{\agfunc_{j_m}}{\xelg'_m}, \ \hact{\helg}{\yelg'})
    \in \Wman
\]
as well, where $(j_1,\ldots,j_m)$ is a permutation of indices $(1,\ldots,m)$ induces by left shift of $\Hgrp$ by itself.
More precisely, $\helg_{j_k} = \helg_k\helg$ for all $k=1,\ldots,m$.

We need to show that $\agfunc_i=\unit{\Ggrp}$ for all $i=1,\ldots,m$, and $\helg = \unit{\Hgrp}$.

Since $\yelg', \hact{\helg}{\yelg'}\in \Vman$, and $\Vman$ is an $\Hgrp$-wandering neighborhood of $\yelg'$, it follows that $\helg=\unit{\Hgrp}$.
Hence,
\[
   \ghact{(\agfunc_1,\ldots,\agfunc_m,\helg)}{(\xelg'_1,\ldots,\xelg'_m,\yelg')} =
   (\gact{\agfunc_{1}}{\xelg'_1}, \ldots, \gact{\agfunc_{m}}{\xelg'_m}, \yelg').
\]
But then $\xelg'_i, \gact{\agfunc_{i}}{\xelg'_i} \in \Uman_i$ for each $i=1,\ldots,m$.
Since $\Uman_i$ is a $\Ggrp$-wandering neighborhood of $\xelg_i$, we obtain that $\agfunc_i=\unit{\Ggrp}$ as well, and thus $\sigma$ is the unit of $\GwrH$.
\end{proof}

\subsection{Main result}
Consider the short exact sequences of $(\Xman,\Gxorb,\xelg)$ and $(\Yman,\Hyorb,\yelg)$:
\begin{align*}
    \pi_1(\Xman,\xelg) \xmonoArrow{\jhomX} \pi_1(\Xman,\Gxorb,\xelg) \xepiArrow{\ehomX} \Ggrp, &&
    \pi_1(\Yman,\yelg) \xmonoArrow{\jhomY} \pi_1(\Yman,\Hyorb,\yelg) \xepiArrow{\ehomY} \Hgrp.
\end{align*}
Then we have the following $(3\times 3)$-diagram whose rows and columns are exact:
\begin{equation}\label{equ:3x3_main}
\begin{aligned}
    \xymatrix@R=2em{
        \bigl(\pi_1\Xman \bigr)^m                       \ar@{^(->}[r]^-{\jhomX^m}  \ar@{^(->}[d] &
        \bigl(\pi_1(\Xman,\Gxorb) \bigr)^m              \ar@{->>}[rr]^-{(\ehomX)^m}     \ar@{^(->}[d] &&
        \Ggrp^m                                         \ar@{^(->}[d] \\
        \bigl( \pi_1\Xman \bigr)^m \times  \pi_1\Yman   \ar@{^(->}[r]  \ar@{->>}[d] &
        \pi_1(\Xman,\Gxorb) \wrm{\ehomY} \pi_1(\Yman,\Hyorb)
                                                        \ar@{->>}[rr]^-{\ehomX\,\wr\,\ehomY}   \ar@{->>}[d]^{j} &&
        \GwrH                                           \ar@{->>}[d] \\
        \pi_1\Yman                                      \ar@{^(->}[r]^-{\jhomY} &
        \pi_1(\Yman,\Hyorb)                             \ar@{->>}[rr]^-{\ehomY}  &&
        \Hgrp
    }
\end{aligned}
\end{equation}
where we omitted base points,
\begin{align*}
    & (\ehomX\,\wr\,\ehomY)(a_1,\ldots,a_m,b)=\bigl(\ehomX(a_1),\ldots,\ehomX(a_m), \ehomY(b) \bigr), &
    & j(a_1,\ldots,a_m,b) = \jhomY(b).
\end{align*}

\begin{subtheorem}\label{th:pi1WGwrH}\sl
Suppose that the action of $\Ggrp$ on $\Xman$ is WD at $\xelg$, and the action of $\Hgrp$ on $\Yman$ is WD at $\yelg$, so $\GwrH$ is WD at $\welg$, see~\eqref{equ:point_w} and Lemma~\ref{lm:GwrH_W_PD}.
Then the short exact sequence of $(\Wman,\GHworb,\welg)$:
\[
    \pi_1(\Wman,\welg) \xmonoArrow{} \pi_1(\Wman,\GHworb,\welg) \xepiArrow{} \GwrH
\]
is isomorphic to the middle horizontal sequence of~\eqref{equ:3x3_main}.
In particular, there is an isomorphism
\begin{equation}\label{equ:pi1W__pi1X_wr_pi1Y__rel}
    \pi_1(\Wman,\GHworb,\welg) \,\cong\, \pi_1(\Xman,\Gxorb,\xelg) \wrm{\dhomY} \pi_1(\Yman,\Hyorb,\yelg).
\end{equation}
\end{subtheorem}
\begin{proof}
Denote
\begin{align*}
    \Xsp&:=\PXxG, &
    \Ysp&:=\PYyH, &
    \Wsp&:=\PWwGH.
\end{align*}
By Corollary~\ref{cor:pi0_PathsAndLoops} we have natural isomorphisms:
\begin{align*}
    \pi_0\Xsp &\cong\pi_1(\Xman,\Gxorb,\xelg), &
    \pi_0\Ysp &\cong\pi_1(\Yman,\Hyorb,\yelg), &
    \pi_0\Wsp &\cong\pi_1(\Wman,\GHworb,\welg).
\end{align*}
It will be more convenient to work with these $\pi_0$-groups.
In particular, the diagram~\eqref{equ:3x3_main} can be written in terms of them, and we thus need to establish the following isomorphism
\[
    \pi_0\Wsp \ \cong \ \pi_0\Xsp\wrm{\ehomY\colon\,\pi_0\Ysp\to\Hgrp}\pi_0\Ysp :=
    \Maps{\Hgrp}{\pi_0\Xsp}\mathop{\rtimes}\limits_{\ehomY}\pi_0\Ysp.
\]

\begin{enumerate}[label={\rm\Alph*)}, ref={\rm\Alph*}, wide, itemsep=2ex, topsep=0.5ex]
\item\label{enum:AA}
We will show in~\ref{enum:AA:1} and~\ref{enum:AA:2} below that {\em there exists a homeomorphism
\[
\zxmyIso: \Cont{I}{\Wman} \ \cong \ \Maps{\Hgrp}{C(I,\Xman)} \times C(I,\Yman),
\]
such that $\zxmyIso(\Wsp) = \Maps{\Hgrp}{\Xsp}\times \Ysp$.}
Hence, $\zxmyIso$ induces a bijection
\begin{equation}\label{equ:eta0_is_iso_1}
    \zxmyIsoZ:\pi_0\Wsp \ \cong \ \pi_0\Maps{\Hgrp}{\Xsp}\times\pi_0\Ysp \ \equiv \ \Maps{\Hgrp}{\pi_0\Xsp}\times\pi_0\Ysp
\end{equation}
of the corresponding sets of path components.
Notice that the latter set has the structure of the group $\pi_0\Xsp\wrm{\ehomY}\pi_0\Ysp$.
We will then show in~\ref{enum:BB:1} and~\ref{enum:BB:2} that~\eqref{equ:eta0_is_iso_1} is an isomorphism of groups:
\begin{equation}\label{equ:eta0_iso}
    \zxmyIsoZ:\pi_0\Wsp \ \cong \ \pi_0\Xsp\wrm{\ehomY}\pi_0\Ysp.
\end{equation}

\begin{enumerate}[label={\rm\ref{enum:AA}\arabic*)}, wide, itemsep=1ex, topsep=0.5ex]
\item\label{enum:AA:1}

Since $I$ and $\Hgrp$ are compact and Hausdorff, and $\Hgrp$ is also discrete (so any map from $\Hgrp$ is continuous), we have (by the exponential law, see Lemma~\ref{cor:pi0_PathsAndLoops}) the canonical homeomorphisms with respect to the corresponding compact open topologies:
\begin{equation}\label{equ:CIHX_HCIX}
\begin{aligned}
    \Cont{I}{\Maps{\Hgrp}{\Xman}} &\equiv
    \Cont{I}{\Cont{\Hgrp}{\Xman}}  \cong
    \Cont{I\times \Hgrp}{\Xman} \cong \\  &\cong
    \Cont{\Hgrp\times I}{\Xman} \equiv
    \Cont{\Hgrp}{\Cont{I}{\Xman}} \equiv
    \Maps{\Hgrp}{\Cont{I}{\Xman}}.
\end{aligned}
\end{equation}
There is also another canonical identification:
\begin{equation}\label{equ:CIW_CIHX_x_CIY}
    \Cont{I}{\Wman} \ \cong \
    \Cont{I}{\Maps{\Hgrp}{\Xman}\times\Yman} \ \cong \
    \Cont{I}{\Maps{\Hgrp}{\Xman}} \times \Cont{I}{\Yman}
\end{equation}
associating to each path $\wpath=(\xpath,\ypath)\colon I \to \Wman=\Maps{\Hgrp}{\Xman}\times\Yman$ the pair of its coordinate functions $(\xpath,\ypath)$.
Hence, we get a homeomorphism
\[
\zxmyIso:
    \Cont{I}{\Wman}
    \ \stackrel{\eqref{equ:CIW_CIHX_x_CIY}}{\cong} \
    \Cont{I}{\Maps{\Hgrp}{\Xman}}\times\Cont{I}{\Yman}
    \ \stackrel{\eqref{equ:CIHX_HCIX}}{\cong} \
    \Maps{\Hgrp}{C(I,\Xman)} \times C(I,\Yman),
\]
defined as follows.
Let $\wpath=(\agpath,\ahpath)\colon I \to \Wman = \Maps{\Hgrp}{\Xman}\times\Yman$ be a path in $\Wman$.
Define the map $\pagpath\colon \Hgrp \to C(I,\Xman)$ by $\pagpath(\helg)(t) := \agpath(t)(\helg)$.
Then
\[ \zxmyIso(\wpath) = (\pagpath, \ahpath).\]

\item\label{enum:AA:2}
We claim that $\zxmyIso(\Wsp) = \Maps{\Hgrp}{\Xsp}\times \Ysp$.
Indeed, let $\wpath=(\agpath,\ahpath) \in C(I,\Wman)$ and $\eta(\wpath)=(\pagpath,\ahpath)$.
Then $\wpath\in\Wsp=\PWwGH$, means that
\begin{equation}\label{equ:w_in_W}
\begin{aligned}
\wpath(0)&=(\agpath(0),\ahpath(0)) = \welg = (\cpathx, \yelg), \\
\wpath(1)&=(\agpath(1),\ahpath(1)) = \ghact{(\agfunc,\helg)}{(\cpathx, \yelg)} =
\bigl( \gprod{\agfunc}{(\cpathx\circ\RightShift{\helg})}, \hact{\helg}{\yelg} \bigr) =
\bigl( \gprod{\agfunc}{\cpathx}, \hact{\helg}{\yelg} \bigr),
\end{aligned}
\end{equation}
for some $(\agfunc,\helg) = \ehomW(\wpath) \in\GwrH$.
In turn, \eqref{equ:w_in_W} is equivalent to the assumption that for each $\kelg\in\Hgrp$ we have that
\begin{equation}\label{equ:ehomW}
\begin{aligned}
    \pagpath(\kelg)(0) &= \agpath(0)(\kelg)=\xelg, & \qquad\qquad
    \ahpath(0) &= \yelg, \\
    \pagpath(\kelg)(1) &= \agpath(1)(\kelg)=\gact{\agfunc(\kelg)}{\xelg} \in \Gxorb, &
    \ahpath(1) &= \hact{\helg}{\yelg} \in \Hyorb,
\end{aligned}
\end{equation}
i.e.\ $\pagpath\in\Maps{\Hgrp}{\Xsp}$ and $\ahpath\in\Ysp$.
Hence, $\zxmyIso(\wpath)=(\pagpath,\ahpath) \in \Maps{\Hgrp}{\Xsp}\times\Ysp$.
\end{enumerate}

\item\label{enum:WW}
We need to prove that~\eqref{equ:eta0_iso} is an isomorphism of groups.
Since it is a bijection, it suffices to check that $\zxmyIsoZ$ is a homomorphism of groups.

\begin{enumerate}[label={\rm\ref{enum:WW}\arabic*)}, wide, itemsep=1ex, topsep=0.5ex]
\item\label{enum:BB:1}
First let us write down explicit formulas for the multiplication in $\pi_0\Wsp$.
Let
\[ \aghpath=(\agpath,\ahpath), \bghpath=(\bgpath,\bhpath)\colon I \to \Wman=\Maps{\Hgrp}{\Xman}\times\Yman
\]
be two paths belonging to $\PWwGH$, and
\begin{align*}
    \cghpath=(\cgpath,\chpath) = \pprod{\aghpath}{\bghpath}
    \ \stackrel{\eqref{equ:prod_xpath_xpath1}}{=} \
    \pconc{\aghpath}{(\gpact{\ehomW(\aghpath)}{\bghpath})}
\end{align*}
be their product.
We need to find precise formulas for the corresponding coordinate functions $\cgpath\colon I\to\Maps{\Hgrp}{\Xman}$ and $\chpath\colon I\to\Yman$.
Since the action of $\GwrH$ is WD at $\welg$, there exists a unique map $\agfunc\colon\Hgrp\to\Xman$ and a unique $\helg\in\Hgrp$ such that
\begin{align*}
    \agpath(1)(\kelg) &= \gact{\agfunc(\kelg)}{\xelg}, \ \kelg\in\Hgrp, &
    \ahpath(1)        &= \hact{\helg}{\yelg},
\end{align*}
so $\ehomW(\aghpath) = (\agfunc, \helg)$, and
$(\gpact{\ehomW(\aghpath)}{\bghpath})(t) =
\ghact{(\agfunc, \helg)}{(\bgpath(t),\bhpath)}=
\bigl(\gprod{\agfunc}{\bgpath(t)\mcirc\RightShift{\helg}}, \hact{\helg}{\yelg}\bigr)$ for all $t\in I$.
Hence,
\begin{equation}\label{equ:prod_in_pi0W}
\begin{aligned}
\cgpath(t)(\kelg) &=
\begin{cases}
\agpath(2t)(\kelg), & t\in[0;\tfrac{1}{2}], \\
\gprod{\phi(\kelg)}{\bgpath(2t-1)(\kelg\helg)}, & t\in[\tfrac{1}{2};1],
\end{cases}
\\
\chpath(t) &=
\begin{cases}
\ahpath(2t), & t\in[0;\tfrac{1}{2}], \\
\hprod{\helg}{\bhpath(2t-1)}, & t\in[\tfrac{1}{2};1].
\end{cases}
\end{aligned}
\end{equation}

\item\label{enum:BB:2}
Now let us describe the multiplication in
\[ \pi_0\Xsp\wrm{\ehomY}\pi_0\Ysp = \Maps{\Hgrp}{\pi_0\Xsp} \rtimes_{\ehomY} \pi_0\Ysp. \]
By definition, each element of that group is a pair $(\agpmap,\ahpmap)$, where $\agpmap\colon\Hgrp\to\pi_0\Xsp$ is some map and $\ahpmap\in\pi_0\Ysp$.
It is thus can be represented by a pair
\[ (\agmap,\ahmap)\in \Maps{\Hgrp}{\Xsp}\times\Ysp = \zxmyIso(\Wsp),\]
such that $\agpmap(\kelg) = [\agmap(\kelg)]$ for all $\kelg\in\Hgrp$, and $\ahpmap = [\ahmap]$, where square brackets mean homotopy classes of the corresponding paths.

Now let $(\agpmap,\ahpmap), (\bgpmap,\bhpmap) \in \pi_0\Xsp\wrm{\ehomY}\pi_0\Ysp$ be two elements, and
\[
    (\cgpmap,\chpmap) := (\agpmap,\ahpmap) (\bgpmap,\bhpmap) =
    \bigl(
        \gprod{\agpmap}{\bgpmap \circ\RightShift{\dhomY(\ahpmap)}},
        \gprod{\ahpmap}{\bhpmap}
    \bigr)
\]
be their product.
Recall that the multiplication here is the point-wise (in $\Hgrp$) concatenation of paths defined by~\eqref{equ:prod_xpath_xpath1}.
Hence, if $(\agmap,\ahmap), (\bgmap,\bhmap) \in \Maps{\Hgrp}{\Xsp}\times\Ysp$ are representatives of those elements, then $(\cgpmap,\chpmap)$ is represented by $(\cgmap,\chmap)\in\Maps{\Hgrp}{\Xsp}\times\Ysp$, where
\begin{align*}
    \cgmap(\kelg) &=
    \pconc{\agmap(\kelg)}{(\gpact{\ehomXx(\agmap(\kelg))}{\bgmap(\kelg)})}, &
    \chmap(\kelg) &=
    \pconc{\ahmap}{(\gpact{\ehomYy(\ahmap)}{\bhmap})}.
\end{align*}
Define $\phi\colon \Hgrp\to\Xman$ by $\phi(\kelg) = \agmap(\kelg)(1)$, and let $\helg = \ahmap(1)$.
Then,
\begin{equation}\label{equ:prod_in_pi0Xwrpi0Y}
\begin{aligned}
    \cgmap(\kelg)(t) &=
\begin{cases}
    \agmap(\kelg)(2t), & t\in[0;\tfrac{1}{2}], \\
    \gprod{\phi(\kelg)}{\bgmap(\kelg\helg)(2t-1)}, & t\in[\tfrac{1}{2};1],
\end{cases}
\\
\chpath(t) &=
\begin{cases}
\ahpath(2t), & t\in[0;\tfrac{1}{2}], \\
\hprod{\helg}{\bhpath(2t-1)}, & t\in[\tfrac{1}{2};1].
\end{cases}
\end{aligned}
\end{equation}

Comparing~\eqref{equ:prod_in_pi0Xwrpi0Y} and~\eqref{equ:prod_in_pi0W} we see that if
\[
    (\agpath,\ahpath), \ \
    (\bgpath,\bhpath), \ \
    (\cgpath,\chpath)=\pprod{(\agpath,\ahpath)}{(\bgpath,\bhpath)} \in \Wsp,
\]
and
\[
    (\agmap,\ahmap) = \zxmyIso(\agpath,\ahpath), \qquad
    (\bgmap,\bhmap) = \zxmyIso(\bgpath,\bhpath), \qquad
    (\cgmap,\chmap) = \pprod{(\agmap,\ahmap)}{(\bgmap,\bhmap)}
\]
then $\zxmyIso(\cgpath,\chpath) = (\cgmap,\chmap)$.
This implies that $\zxmyIsoZ$ is a homomorphism, and therefore an isomorphism.
\qedhere
\end{enumerate}
\end{enumerate}
\end{proof}

\begin{subcorollary}\label{cor:coverings}\sl
Suppose $p\colon\Xman\to\Xman/\Ggrp$ and $q\colon\Yman\to\Yman/\Hgrp$ are regular covering maps, and $\Hgrp$ is finite of order $m$.
Then $r\colon \Wman\to\Wman/(\GwrH)$ is also a regular covering map, and we have an isomorphism:
\begin{equation}\label{equ:pi1W__pi1X_wr_pi1Y}
\pi_1\bigl(\Wman/(\GwrH)\bigr) \ \cong \ \pi_1 (\Xman/\Ggrp) \, \wrm{\dhomY} \, \pi_1(\Yman/\Hgrp).
\end{equation}
In particular, if $p$ and $q$ are universal covering maps, then $\pi_1\bigl(\Wman/(\GwrH)\bigr) \cong \Ggrp \wr \Hgrp$.
\end{subcorollary}
\begin{proof}
Let $\xelg\in\Xman$ and $\yelg\in\Yman$ be any points, and $\welg\in(\underbrace{\xelg,\ldots,\xelg}_{m},y) \in \Wman = \Xman^m\times\Yman$.
As $p,q,r$ are regular coverings, the corresponding actions of $\Ggrp$, $\Hgrp$ and $\GwrH$ are WD at $\xelg$, $\yelg$ and $\welg$ respectively.
In particular, and we have isomorphisms $\pi_1(\Xman/\Ggrp)\cong\pi_1(\Xman,\Gxorb)$, $\pi_1(\Yman/\Hgrp)\cong\pi_1(\Yman,\Hyorb)$, and $\pi_1\bigl(\Wman/(\GwrH)\bigr)\cong\pi_1(\Wman,\GHworb)$, where $x,y,w$ are any points as in Theorem~\ref{th:pi1WGwrH}.
Now by Theorem~\ref{th:pi1WGwrH}, the isomorphism~\eqref{equ:pi1W__pi1X_wr_pi1Y__rel} gives an isomorphism~\eqref{equ:pi1W__pi1X_wr_pi1Y}.
\end{proof}

Let $\Tor{0}$ be a singleton, and $\Tor{n} =(\Circle)^n$, $n\geq1$, be the $n$-torus.

\begin{subexample}\rm
Let $p\colon\Xman\to\Xman/\Ggrp$ be a covering map and $q\colon\Circle\to\Circle$, $q(\welg)=\welg^m$, be the quotient map corresponding to the action of $\Hgrp=\bZ_m$ on $\Yman=\Circle$ by rotations.
Then the boundary homomorphism $\dhomY\colon\pi_1\Circle=\bZ\to\bZ_m$ of $q$ is given by $\dhomY(a)=a\bmod m$, so the short exact sequence of $q$ is $\seqC{m}\colon m\bZ\monoArrow\bZ\epiArrow\bZ_{m}$.

In particular, due to Example~\ref{enum:W6}, $\pi_1(\QXG)\,\wrm{\dhomY}\,\pi_1(\QYH)$ is the same as $\pi_1(\QXG)\,\wrm{m}\,\bZ$.
Moreover, $\Wman = \Xman^{m}\times\Circle$, and the action $\Ggrp\wr\bZ_m$ on $\Wman$ is given by
\[
(\gelg_0,\ldots,\gelg_{m-1}, \kelg) (x_0,\ldots,x_{m-1}, \welg) =
\bigl(
\gelg_0(x_{\kelg}),\gelg_1(x_{1+\kelg}),\ldots,\gelg_{m-1}(x_{\kelg-1}), \welg e^{2\pi i \kelg/m}\bigr),
\]
where the lower indices in $\xelg_{i}$ are taken modulo $m$.
Hence, by Theorem~\ref{th:pi1WGwrH}, the short exact sequence of the covering map $r\colon \Wman\to\Wman/(\GwrH)$ is
\[
    \seqWrm{\pSeq}{m}\colon
    (\pi_1\Xman)^m\times m\bZ   \xmonoArrow{~p^m\times q~}
    \pi_1(\QXG)\,\wrm{m}\,\bZ   \xepiArrow{~\dhomX\,\wr\,\dhomY~}
    \Ggrp\wr\bZ_m,
\]
so in particular $\pi_1\Wman \cong \pi_1(\Xman/\Ggrp) \wrm{m} \bZ$.
\end{subexample}

\begin{subexample}\rm
More generally, let $m,n\geq1$ and $q\colon\Tor{2} \to \Tor{2}$ be the covering map given by $q(w_1,w_2)=(w_1^m,w_2^n)$.
It corresponds to the free action of $\Hgrp=\bZ_{m}\times\bZ_{n}$ on $\Yman=\Tor{2}$ defined by $(a,b)(w_1,w_2)=(w_1 e^{2\pi i a/m}, w_2 e^{2\pi i b/n})$, $a,b\in\bZ_{m}\times\bZ_{n}$.
Then the short exact sequence of $q$ is
\[
    \seqC{m}\times\seqC{n}:
        m\bZ\times n\bZ \xmonoArrow{~q~}
        \bZ\times\bZ    \xepiArrow{~\dhomY:\,(a,\,b) \,\mapsto\, (a\bmod m,\, b\bmod n)~}
        \bZ_{m} \times \bZ_{n}.
\]
Now, due to Example~\ref{enum:W8}, $\pi_1(\QXG)\,\wrm{\dhomY}\,\pi_1(\QYH)$ is the same as $\pi_1(\QXG)\,\wrm{m,n}\,\bZ^2$.
We also have a PD action of $\Ggrp\wr\Hgrp$ on $\Wman = \Xman^{mn}\times\Yman$ defined by
\begin{multline*}
    \bigl( \{\gelg_{i,\,j}\}_{i\in\bZ_m,\,j\in\bZ_n}, \, a,       \, b\bigr) \cdot
    \bigl( \{\xelg_{i,\,j}\}_{i\in\bZ_m,\,j\in\bZ_n}, \, \welg_1, \, \welg_2\bigr)
        = \\ =
    \bigl(
       \{\gelg_{i,\,j}(x_{i+a,\,j+b})\}_{i\in\bZ_m,\,j\in\bZ_n}, \,
       \welg_1 e^{2\pi i a/m}, \,
       \welg_2 e^{2\pi i b/n}
    \bigr).
\end{multline*}
Then by Theorem~\ref{th:pi1WGwrH}, the short exact sequence of $r\colon \Wman\to\Wman/(\GwrH)$ is
\[
    \seqWrm{\pSeq}{m,n}\colon
    (\pi_1\Xman)^m\times m\bZ \times n\bZ  \xmonoArrow{~p^m\times q~}
    \pi_1(\QXG)\,\wrm{m,n}\,\bZ^2          \xepiArrow{~\dhomX\,\wr\,\dhomY~}
    \Ggrp\wr(\bZ_m\times\bZ_n).
\]
In particular, $\pi_1\Wman \cong \pi_1(\Xman/\Ggrp) \wrm{m,n} \bZ^2$.
\end{subexample}

\begin{subexample}\label{exmp:action:wrn}\rm
Similarly, for a product of finitely many cyclic groups
\[
    \Hgrp=\bZ_{m_1}\times\cdots\times\bZ_{m_n}
\]
with $m=m_1\cdots m_n$, one can similarly define a free action of $\GwrH$ on $\Wman=\Xman^{m}\times\Tor{n}$ such that the short exact sequence of the covering map $r\colon \Wman\to\Wman/(\GwrH)$ is
\[
    \seqWrm{\pSeq}{m_1,\ldots,m_k}\colon
    (\pi_1\Xman)^m\times \prod\limits_{i=1}^{k} m_i\bZ  \xmonoArrow{~~~~}
    \pi_1(\QXG)\,\wrm{m_1,\ldots,m_k}\,\bZ^k     \xepiArrow{~~~~}
    \Ggrp\wr\prod\limits_{i=1}^{k} \bZ_{m_i},
\]
so in particular $\pi_1(\Wman/(\GwrH)) \ \cong \ \pi_1(\QXG)\wrm{m_1,\ldots,m_n}(\bZ_{m_1}\times\cdots\times\bZ_{m_n})$.
\end{subexample}

\section{Classes of short exact sequences}\label{sect:classes_B}
Iterated wreath products naturally act on trees, e.g.~\cite{Boudec:ETDS:2021}.
That observation was made in 1869 by C.~Jordan~\cite{Jordan:JRAM:1869} who probably introduced wreath products, see~\cite[p.~209]{Polya:AM:1937}.
Let $\mathcal{T}$ be the set of isomorphism classes of groups of automorphisms of finite trees.
Jordan proved that every group from $\mathcal{T}$ is obtained from the unit group by finitely many operations of direct product and wreath product of the form $\cdot\wrm{X_n} S_n$, where $X_n=\{1,\ldots,n\}$ and $S_n$ is the permutation group of $X_n$.
We will now define in a similar way two classes of groups generated by wreath products $\cdot\wr_{m_1,\ldots,m_k}\bZ^k$ and $\cdot\wr\prod_{i=1}^{k}\bZ_{m_i}$ and related by certain short exact sequences.

\begin{definition}\label{def:seq:classBn}\sl
For $n\geq1$ let $\classB_{n}$ be the minimal set of isomorphism classes of short exact sequences having the following properties (see Section~\ref{sect:short_ex_sequence} for notations):
\begin{enumerate}[label={\rm(\alph*)}]
\item\label{enum:seq:classB:1}
the sequence $\seqC{0}:1\monoArrow1\epiArrow1$ belongs to $\classB_{n}$;

\item\label{enum:seq:classB:x}
if $\pSeq,\qSeq\in\classB$, then $\pSeq\times\qSeq\in\classB_{n}$;

\item\label{enum:seq:classB:wr}
if $\pSeq\in\classB$ and $m_1,\ldots,m_k\geq1$ is a finite collection of natural numbers with $k\leq n$, then $\seqWrm{\pSeq}{m_1,\ldots,m_k} \ \in \classB_{n}$.
\end{enumerate}
\end{definition}

Due to~\ref{enum:seq:classB:1}, the set $\classB_{n}$ is non-empty, and its minimality means that every short exact sequence $\pSeq\in\classB_{n}$ is obtained from $\seqC{0}$ by finitely many operations of direct product and wreath product of the form $\seqWrm{\cdot}{m_1,\ldots,m_k}$ with $k\leq n$.
For example,
\begin{align*}
\seqC{k} \in \classB_{1}, (k\geq1), \quad
\seqWrm{(\seqC{2}\times\seqC{9})}{5,7}\in\classB_{2}, \quad
\seqWrm{\bigl(\seqC{32}\times(\seqWrm{\seqC{2}}{5,7,3})\bigr) }{9,34,6,2} \in\classB_{4}.
\end{align*}
It is also evident that $\classB_{k} \subset \classB_{l}$ for $k<l$.
Put $\classB := \cup_{n\geq1}\classB_{n}$.

\begin{remark}\rm
Let $\classBB$ and $\classBC$ be the sets of isomorphism classes of middle and right groups appearing in the short exact sequences $\pSeq\colon A \monoArrow B \epiArrow C \in \classB$.
Then it follows from the definition that $\classBB$ (resp.\ $\classBC$) is the minimal set of classes containing the unit group $1$ and closed under direct product and wreath products $\cdot\wrm{m_1,\ldots,m_k}\bZ^k$ (resp.\ $\cdot\wr \prod\limits_{i=1}^{k}\bZ_{m_i}$).
Thus, $\classBB$ and $\classBC$ are defined similarly to $\mathcal{T}$.
\end{remark}

Any short exact sequence isomorphic to a sequence of the form $\pSeq\colon\bZ^n \monoArrow B \epiArrow C$ with finite $C$ and some $n\geq0$ will be called \term{crystallographic}.
In this case the middle group $B$ will also be called \term{crystallographic}.
We will also say that $\pSeq$ as well as $B$ are \term{Bieberbach}, whenever $B$ is torsion free.
In~\cite[Lemmas~2.2 \& 2.6]{Maksymenko:TA:2020} the author proved that every sequence from $\classB_{1}$ is Bieberbach.

\begin{remark}\label{rem:what_to_fix}\rm
In~\cite[Theorem 2.5]{Maksymenko:TA:2020} the following statement is mentioned as \term{Bieberbach theorem} and referred to~\cite[Corollary~5.1]{Charlap:BG:1986} for its proof:
\term{for every Bieberbach sequence $\pSeq: \bZ^n \monoArrow B \epiArrow C$ there exists a free action of $C$ on a torus $\Tor{n}$ such that the exact sequence of that covering $\bZ^{n} \monoArrow  \pi_1(\Tor{n}/C) \epiArrow C$ is isomorphic to $\pSeq$ relatively to the last term $C$}.
In particular, $\pi_1(\Tor{n}/C) \cong B$.
However, Bieberbach theorem claims existence of such actions only under assumption that $\bZ^n$ is a \term{maximal} abelian subgroup.

Let us show that the above statement is true even for non-maximal free abelian subgroups, which will give lacking arguments.
This is essential, since \cite[Theorem~5.10]{Maksymenko:TA:2020} and the orientable part of~\cite[Corollary~1.3]{Maksymenko:TA:2020} are based on the variant of Bieberbach theorem for non-maximal subgroups.
For example, in the sequence $\seqC{m}\colon m\bZ \monoArrow \bZ \epiArrow \bZ_m$ with $m\geq2$, the abelian subgroup $m\bZ$ of finite index is non-maximal, however there is a free action of $\bZ_m$ on $\Tor{1}=\Circle$ (by rotations by $2\pi/m$) such that $\seqC{m}$ is isomorphic with the short exact sequence of the quotient map $p\colon\Tor{1} \to \Tor{1}/\bZ_m = \Tor{1}$.

Thus, let $A \monoArrow B \epiArrow C$ be a Bieberbach sequence, where $A$ is a free abelian subgroup of $B$ of rank $n$.
If $A$ is non-maximal, take any maximal abelian subgroup $A'\subset B$ containing $A$.
In fact, see A.~Vasquez~\cite[Theorem~3.1]{Vasquez:JDG:1970}, $A'$ is the centralizer of the normal subgroup $A$ in $B$, and therefore it is normal as well.
Since $A'/A$ is finite, $A' \cong\bZ^{n}$ and we get another Bieberbach sequence $\qSeq\colon A' \monoArrow B \epiArrow C'$.
Then by the ``maximal subgoup'' case of Bieberbach theorem, there exists an action of $C':=B/A'$ on $T':=\Tor{n}$ such that the exact sequence of the corresponding covering map $q\colon T' \to T'/C'$ is isomorphic to $\qSeq$.
Let $r\colon T = \Tor{n}\to T'$ be the covering map corresponding to the subgroup $A\subset A'$, and $p=q\circ r\colon T \xrightarrow{r} T' \xrightarrow{q} T'/C'$.
Then $p(\pi_1(T)) = A$, so $p$ is a regular covering, its short exact sequence is isomorphic with $\pSeq$, and we can identify $T'/C'=T/C$.
In particular, this gives the required action of $C$ on $T$ such that $\pi_1 (T/C) \cong B$.
\end{remark}

The following lemma provides an explicit description of the actions on tori for the short sequences from class $\classB$.

\begin{lemma}[{cf.~\cite[Lemma~2.6]{Maksymenko:TA:2020}}]\label{lm:seq_ABC}\sl
For each $\pSeq:A\monoArrow B \epiArrow C \in\classB$ the following statements hold.
\begin{enumerate}[label={\rm(\roman*)}, leftmargin=*]
\item\label{enum:seq_ABC:prop}
$\pSeq$ is Bieberbach, so $A\cong \bZ^n$ for some $n$, $B$ is torsion free, and $C$ is finite.
Moreover, $B,C$ are also solvable.

\item\label{enum:seq_ABC:act}
There exists a free action of the group $C$ on $\Tor{n}$ such that $\pSeq$ is isomorphic to the short exact sequence $\pi_1\Tor{n} \monoArrow \pi_1(\Tor{n}/C) \epiArrow C$ of the corresponding covering map $p\colon\Tor{n}\to \Tor{n}/C$.
\end{enumerate}
Therefore, since $\Tor{n}/C$ is aspherical, i.e.\ it is an Eilenberg-MacLane space $K(B,1)$, any other aspherical path connected topological space $\Xman$ with $\pi_1\Xman\cong B$ is weakly homotopy equivalent to $\Tor{n}/C$.
\end{lemma}
\begin{proof}
Let $\classB' \subset \classB$ be the subset of isomorphism classes of sequences having properties~\ref{enum:seq_ABC:prop} and~\ref{enum:seq_ABC:act}.
We need to prove that $\classB'=\classB$.

It suffices to show that $\classB'$ satisfies conditions~\ref{enum:seq:classB:1}-\ref{enum:seq:classB:wr} of Definition~\ref{def:seq:classBn}.
Since $\classB$ is the minimal class with those properties, we will then have that $\classB'=\classB$.

\ref{enum:seq:classB:1}
Let $\pSeq=\seqC{0}:\bZ^{0} \monoArrow 1 \epiArrow 1$.
Then~\ref{enum:seq_ABC:prop} trivially holds for $\seqC{0}$.
Also, the latter group $C=1$ freely acts on the singleton $\Tor{0}$ so that the short exact sequence of the corresponding covering map $p\colon\Tor{0} \to \Tor{0}/C$ is $\seqC{0}$.
This also proves~\ref{enum:seq_ABC:act}, so $\seqC{0} \in \classB'$.

\ref{enum:seq:classB:x}
Let $\pSeq_i: \bZ^{n_i} \monoArrow B_i \epiArrow C_i \in\classB'$, $i=1,2$.
Thus, each $B_i$ is solvable torsion free, $C_i$ is finite solvable and freely acts on $\Tor{n_i}$ so that $\pSeq_i$ is isomorphic with the short exact sequence of the covering map $p_i\colon\Tor{n_i}\to\Tor{n_i}/C_i$.
Then, $B_1\times B_2$ is also solvable torsion free, and the product $C_1\times C_2$ is finite solvable and naturally acts on $\Tor{n_1}\times\Tor{n_2}=\Tor{n_1+n_2}$ so that the short exact sequence of the covering map $p\colon\Tor{n_1+n_2}\to\Tor{n_1+n_2}/(C_1 \times C_2)$ is isomorphic with $\pSeq_1\times\pSeq_2\colon\bZ^{n_1+n_2} \monoArrow B_1\times B_2 \epiArrow C_1\times C_2$.
Thus, $\pSeq_1\times\pSeq_2\in\classB'$.

\ref{enum:seq:classB:wr}
Finally, let $\pSeq\colon\bZ^{n}\monoArrow B \epiArrow C \in\classB'$, $m_1,\ldots,m_k\geq1$ be any natural numbers, and $m=m_1\cdots m_k$.
In particular, we have a free action of the group $C$ on $\Tor{n}$ such that $\pSeq$ is short exact sequence of the covering map $p\colon\Tor{n} \to \Tor{n}/C$.
Then $\Hgrp := C\wr\prod\limits_{i=1}^{k}\bZ_{m_i}$ is finite and solvable as well as $C$.
Moreover, by Example~\ref{exmp:action:wrn}, there exists a free action of $\Hgrp$ on the torus $(\Tor{n})^m\times \Tor{k}=\Tor{nm+k}$ such that the short exact sequence of the covering map $q\colon\Tor{nm+k} \to \Tor{nm+k}/\Hgrp$ is isomorphic with
\[
    \seqWrm{\pSeq}{m_1,\ldots,m_k}:
    \bZ^{nm+k} \xmonoArrow{~~} B\wrm{m_1,\ldots,m_k}\bZ^k \xepiArrow{~~}
    C\wr\prod\limits_{i=1}^{k} \bZ_{m_i}.
\]
Note that $B\wrm{m_1,\ldots,m_k}\bZ^k$ is also solvable as well as $B$.
Moreover, by the arguments similar to~\cite[Lemma~2.2]{Maksymenko:TA:2020} for $k=1$, $B\wrm{m_1,\ldots,m_k}\bZ^k$ is also torsion free.
Hence, $\seqWrm{\pSeq}{m_1,\ldots,m_k}\in\classB$.
\end{proof}

\section{Homotopy types of orbits of smooth functions on surfaces}
\label{sect:orbits}
Let $\Mman$ be a compact surface and $\Pman$ be either the real line $\bR$ or the circle $\Circle$.
For a closed subset $\Xman\subset\Mman$ denote by $\Diff(\Mman,\Xman)$ the group of all smooth ($\Cinfty$) diffeomorphisms of $\Mman$ fixed on $\Xman$.
Then the group $\Diff(\Mman,\Xman)$ acts on the space $\Ci{\Mman}{\Pman}$ by the following rule: if $\dif\in\Diff(\Mman,\Xman)$ and $\func\in \Ci{\Mman}{\Pman}$, then the result of the action of $\dif$ on $\func$ is the composition map $\func\circ\dif:\Mman\to\Pman$.
For $\func\in \Ci{\Mman}{\Pman}$ let
\begin{align*}
\Stabilizer{\func,\Xman} &= \{\dif \in \Diff(\Mman,\Xman) \mid \func \circ \dif = \func \}, &
\Orbit{\func,\Xman} &= \{\func \circ \dif \mid \dif \in \Diff(\Mman,\Xman) \}
\end{align*}
be respectively the \textit{stabilizer} and the \textit{orbit} of $\func$ under that action.
It will be convenient to say that elements of $\Stabilizer{\func,\Xman}$ \term{preserve} $\func$.
Endow the above spaces with the corresponding strong $C^{\infty}$ topologies and denote by $\DiffId(\Mman,\Xman)$ and $\StabilizerId{\func,\Xman}$ the corresponding path components of $\id_{\Mman}$ in $\Diff(\Mman,\Xman)$ and $\Stabilizer{\func,\Xman}$, and by $\OrbffX$ the path component of $\OrbfX$ containing $\func$.
If $\Xman$ is empty, then it will be omitted from notation.

Let also $\Cid{\Mman}{\Pman} \subset \Ci{\Mman}{\Pman}$ be the subset consisting of maps $\func\colon\Mman\to\Pman$ satisfying the following axiom:
\begin{enumerate}[leftmargin=*, topsep=1ex, parsep=1ex, label={\AxBd}]
\item\label{axiom:bd}
\it
$\func$ takes a constant value at every connected component of $\partial\Mman$ and has no critical points on $\partial\Mman$.
\end{enumerate}
A map $\func \in \Cid{\Mman}{\Pman}$ is called \term{Morse} if all its critical points are non-degenerate.
Then the set $\Morse(\Mman,\Pman)$ of Morse maps is open and everywhere dense in $\Cid{\Mman}{\Pman}$.
Further, let $\FSP{\Mman}{\Pman}$ be the subset of $\Cid{\Mman}{\Pman}$ consisting of maps $\func$ satisfying one more axiom:
\begin{enumerate}[leftmargin=*, topsep=1ex, parsep=1ex, label={\AxCrPt}]
\item\label{axiom:sing}
\it
for every critical point $z$ of $\func$, there are local coordinates in which $\func$ is a homogeneous polynomial $\bR^2\to\bR$ of degree $\geq2$ without multiple factors.
\end{enumerate}
By Morse Lemma every non-degenerate singularity is $\Cinfty$ equivalent to a homogeneous polynomial without multiple factors $\pm x^2\pm y^2$ and thus satisfies~\ref{axiom:sing}.
This means that $\Morse(\Mman,\Pman) \subset \FSP{\Mman}{\Pman}$.
Notice that~\ref{axiom:sing} also implies that each critical point of $\func\in\FSP{\Mman}{\Pman}$ is isolated, whence the set of critical points of $\func$ is finite.
By~\cite{ChurchTimourian:PJM:1973, Dancer:2:JRAM:1987, Prishlyak:TA:2002}, for every isolated critical point $z$ of a $\Cr{3}$ function $\func\colon\bR^2\to\bR$ the local \term{topological structure} of level-sets of $\func$ near $z$ is realized by level sets of homogeneous polynomial without multiple factors.
Thus, $\FSP{\Mman}{\Pman}$ contains not only ``all typical'' (i.e.\ Morse) maps, but also maps with all possible topological types of critical points.

Let $\func\in\FSP{\Mman}{\Pman}$.
A connected component $\Kman$ of a level-set $f^{-1}(c)$, $c\in\Pman$, will be called a \term{contour} (of $\func$).
We also call $\Kman$ \term{regular} if it contains no critical points, and \term{critical} otherwise.
Then a compact submanifold $\Xman\subset\Mman$ whose connected components have dimensions $1$ and $2$ will be said \term{$\func$-saturated} if it is a union of contours of $\func$.
In particular, every regular contour of $\func$ and $\partial\Mman$ are $\func$-saturated.
Also, if $a,b\in\Pman$ are regular values of $\func$, and $J \subset \Pman$ an interval with $\partial J = \{a,b\}$, then $\func^{-1}(J)$ is $\func$-saturated.

The study of homotopy types of stabilizers and orbits of Morse maps $\func\in\Morse(\Mman,\Pman)$ was initiated in~\cite{Maksymenko:AGAG:2006}.
In a series of papers~\cite{Maksymenko:AGAG:2006, Maksymenko:TrMath:2008, Maksymenko:ProcIM:ENG:2010, Maksymenko:UMZ:ENG:2012,Maksymenko:TA:2020} it was proved that for $\func\in\FSP{\Mman}{\Pman}$ and an $\func$-saturated submanifold $\Xman$ the following statements hold.
\begin{enumerate}[leftmargin=*,label={\rm(\arabic*)}, itemsep=0.5ex]
\item
The map $p\colon\Diff(\Mman,\Xman)\to\Orbit{\func,\Xman}$ is a Serre fibration with fiber $\Stabilizer{\func}$, \cite[Th.~2.1]{Maksymenko:AGAG:2006}, \cite[Th.~5.1]{Maksymenko:ProcIM:ENG:2010}.

\item
The path component $\StabilizerId{\func}$ is contractible except for few types of maps for which it is homotopy equivalent to the circle, \cite[Th.~1.3]{Maksymenko:AGAG:2006},
\cite[Th.~5.1]{Maksymenko:ProcIM:ENG:2010}, \cite[Th.~2.1]{Maksymenko:UMZ:ENG:2012}.

\item\label{enum:BS:bieb_seq}
There is a short exact sequence $\pi_1\DiffId(\Mman,\Xman)\times\bZ^{k} \monoArrow \pi_1\OrbitPathComp{\func,\Xman}{\func} \epiArrow \Ggrp$, where $\Ggrp$ is a finite group and $k\geq0$ both depending on $\func$ and $\Xman$.

\item
$\OrbitPathComp{\func}{\func}$ is aspherical iff either $\Mman\not=\Sphere,\PrjPlane$ or $\Xman\not=\varnothing$, \cite[Th.~1.5]{Maksymenko:AGAG:2006},  \cite[Th.~2.3]{Maksymenko:UMZ:ENG:2012}.
In that case
\begin{itemize}[leftmargin=*]
\item $\pi_1\DiffId(\Mman,\Xman)$ is a free abelian of rank $\leq 2$.
In fact, $\pi_1\DiffId(\Mman,\Xman)=0$ iff either $\Xman\not=\varnothing$ or $\chi(\Mman)<0$;
$\pi_1\DiffId(\Tor{2})=\bZ^2$; and in all other cases $\pi_1\DiffId(\Mman)=\bZ$.
\item
$\pi_1\OrbitPathComp{\func,\Xman}{\func}$ is torsion free, \cite[Lemma~2.2]{Maksymenko:TA:2020}.
\end{itemize}
Hence, the short exact sequence from~\ref{enum:BS:bieb_seq} is Bieberbach.

\item\label{enum:BS:generic_morse}
If $\func$ is a \term{generic} Morse map (i.e.\ it takes distinct values at distinct critical points), then $\Ggrp$ is trivial and expect for few cases $\OrbitPathComp{\func}{\func}$ has the homotopy type of $\Tor{k} \times R$, where $k\geq0$ and $R=SO(3)$ if $M=S^2,\PrjPlane$ and $R$ is a point otherwise, \cite[Th.~1.5]{Maksymenko:AGAG:2006}.
\end{enumerate}

E.~Kudryavtseva~\cite{KudryavtsevaPermyakov:MatSb:2010, Kudryavtseva:MathNotes:2012, Kudryavtseva:SpecMF:VMU:2012, Kudryavtseva:MatSb:2013, Kudryavtseva:ENG:DAN2016} studied the homotopy type of the space of Morse functions on compact surfaces and rediscovering ideas from~\cite{Maksymenko:AGAG:2006} extended the above results, see~\cite[footnote after Theorem~1.2]{Maksymenko:TA:2020}.
She proved that if $\Mman$ is orientable, then under additional mild assumptions on $\func\in\Morse(\Mman,\bR)$ there exists a free action of $\Ggrp$ on $\Tor{n}$ such that $\OrbitPathComp{\func}{\func}$ is homotopy equivalent to $(\Tor{n}/\Ggrp)\times R$, where $R$ is the same as above, which generalizes~\ref{enum:BS:generic_morse}.

Further in~\cite[Section~5]{Maksymenko:TA:2020} the author described a precise algebraic structure of the Bieberbach sequence~\ref{enum:BS:bieb_seq}, which allowed to explicitly compute all its groups.
Also, in a series of papers with B.~Feshchenko that sequence~\ref{enum:BS:bieb_seq} was also computed for $\Mman=\Tor{2}$.
Those results can be formulated as follows:
\begin{theorem}[{\cite{Maksymenko:TA:2020, MaksymenkoFeshchenko:MS:2015, MaksymenkoFeshchenko:MFAT:2015, Feshchenko:Zb:2015, KravchenkoFeshchenko:MFAT:2020, Feshchenko:PIGC:2021}}]\label{th:fix_O_TG}\sl
Let $\Mman$ be a compact orientable surface, $\func\in\FSP{\Mman}{\Pman}$, $\Xman\subset\Mman$ an $\func$-saturated submanifold, and
\begin{equation}\label{eu:bieb_fX}
\pSeq\colon \bZ^{n}\monoArrow\pi_1\OrbitPathComp{\func,\Xman}{\func}\epiArrow\Ggrp
\end{equation}
the Bieberbach sequence of $(\func,\Xman)$, see~\ref{enum:BS:bieb_seq}.
\begin{enumerate}[label={\rm(\alph*)}]
\item
If $\Mman$ is distinct from $\Tor{2}$ and $\Sphere$, then~\eqref{eu:bieb_fX} belongs to $\classB_1$.
\item
If $\Mman=\Tor{2}$, then~\eqref{eu:bieb_fX} belongs to $\classB_2$.
\end{enumerate}
In both cases, a sequence of operations of direct products and wreath products of the form $\seqWrm{\cdot}{m}$ and $\seqWrm{\cdot}{m_1,m_2}$ generating~\eqref{eu:bieb_fX} can be explicitly written down via $\func$ and $\Xman$.
\qed
\end{theorem}
The following statement extends~\cite[Theorem~5.10]{Maksymenko:TA:2020} to all orientable surfaces distinct from $S^2$ (the new statement corresponds to $M=\Tor{2}$).
\begin{corollary}[{c.f.\ \cite{Kudryavtseva:MatSb:2013, Maksymenko:TA:2020}}]\sl
\label{cor:Of_Tn_G}
Under notation of Theorem~\ref{th:fix_O_TG} there is a free action of $\Ggrp$ on $\Tor{n}$ such that $\pSeq$ coincides with the short exact sequence of the covering map $p\colon\Tor{n}\to\Tor{n}/\Ggrp$.
In this case we have the following weak homotopy equivalence $\OrbitPathComp{\func,\Xman}{\func}\simeq\Tor{n}/\Ggrp$.
\end{corollary}
\begin{proof}
Since $\pSeq\in\classB$, the required action of $\Ggrp$ on $\Tor{n}$ is guaranteed by Lemma~\ref{lm:seq_ABC}.
In particular, we get an isomorphism $\pi_1\OrbitPathComp{\func,\Xman}{\func}\cong\pi_1(\Tor{n}/\Ggrp)$.
As $\OrbitPathComp{\func,\Xman}{\func}$, $\Tor{n}$, and therefore its quotient $\Tor{n}/\Ggrp$ are aspherical, the latter isomorphism is induced by a some weak homotopy equivalence $\OrbitPathComp{\func,\Xman}{\func}\,\simeq\,\Tor{n}/\Ggrp$.
\end{proof}

\subsection*{Acknowledgments}
The author is grateful for the financial support within the program of support for priority research and technical (experimental) developments of the Section of Mathematics of the NAS of Ukraine for 2022-2023.
Project ``Innovative methods in the theory of differential equations, computational mathematics and mathematical modeling'', No.~7/1/241.

The author thanks D.~Bolotov for fruitful discussions of the relative $\pi_1$ set.

The author is also very thankful to the Referee of this paper for careful reading and a lot of suggestions which allowed to improve the exposition.


\end{document}